\documentclass[11pt]{article}
\usepackage{amsmath,amsfonts,amsthm,amssymb, color}
\usepackage{graphicx}
\allowdisplaybreaks

\def\ds{\displaystyle}
\def\forall{\hbox{for all}~}
\def\L{{\bf L}}

\def\bfr{{\bf r}}

\def\ve{\varepsilon}

\def\A{{\cal A}}
\def\J{{\cal J}}

\def\R{\mathbb{R}} 

\def\vp{\varphi}

\def\TV{\hbox{Tot.Var.}}

\def\v{\vskip 1em}
\def\O{{\cal O}}
\def\D{{\cal D}}
\def\C{{\cal C}}

\def\Q{{\cal Q}}
\def\I{{\cal I}}
\def\J{{\cal J}}

\def\Ups{\Upsilon}
\def\bega{\begin{array}}
\def\enda{\end{array}}
\def\begi{\begin{itemize}}
\def\endi{\end{itemize}}
\def\ov{\overline}

\def\Tilde{\widetilde}

\def\meas{\hbox{meas}}
\def\bel{\begin{equation}\label}
\def\eeq{\end{equation}}
\def\sqr#1#2{\vbox{\hrule height .#2pt
\hbox{\vrule width .#2pt height #1pt \kern #1pt
\vrule width .#2pt}\hrule height .#2pt }}
\def\square{\sqr74}
\def\endproof{\hphantom{MM}\hfill\llap{$\square$}\goodbreak}

\usepackage{mathtools}

\newtheorem{theorem}{Theorem}[section]
\newtheorem{cor}[theorem]{Corollary}
\newtheorem{lemma}[theorem]{Lemma}
\newtheorem{remark}[theorem]{Remark}
\newtheorem{prop}[theorem]{Proposition}

\numberwithin{equation}{section}

\definecolor{darkgreen}{rgb}{0.0, 0.5, 0.0}

\begin{document}
\title{\bf Uniqueness  Domains for $\L^\infty$  Solutions of \\
$2\times 2$ 
Hyperbolic Conservation Laws}

\author{Alberto Bressan$^{(*)}$, Elio Marconi$^{(**)}$,  and
Ganesh Vaidya$^{(*)}$\\~~\\
 {\small $^{(*)}$~Department of Mathematics, Penn State University,}
 \\
 {\small $^{(**)}$~Department of Mathematics 
 ``Tullio Levi-Civita", University of Padova,}\\~~\\
{\small E-mails: axb62@psu.edu,~~elio.marconi@unipd.it, ~gmv5228@psu.edu}
}
\maketitle
\v
\begin{abstract} For a genuinely nonlinear $2\times 2$ hyperbolic system of conservation laws, 
assuming that the initial data have small $\L^\infty$ norm but possibly unbounded total variation,
the existence of global solutions was proved in a classical paper by Glimm and Lax (1970).
In general, the total variation of these solutions decays like $t^{-1}$.
Motivated by the theory of fractional domains for linear analytic semigroups, we consider here
solutions with faster decay rate: $\TV\bigl\{u(t,\cdot)\bigr\}\leq C t^{\alpha-1}$.
For these solutions, a uniqueness theorem is proved.  Indeed, as the initial data range over 
a domain of functions with $\|\bar u\|_{\L^\infty} \leq\ve_1$ small enough, solutions with fast decay
yield a H\"older continuous semigroup.  The H\"older exponent can be taken arbitrarily close to 1
by further shrinking the value of $\ve_1>0$. An auxiliary result identifies a class of initial data
whose solutions have rapidly decaying total variation. 
\end{abstract}

\section{Introduction}
\label{sec:1}
\setcounter{equation}{0}
Consider a strictly hyperbolic system of conservation laws
in one space dimension
\bel{1}u_t+ f(u)_x~=~0,\eeq
with initial data
\bel{2} u(0,x)~=~\ov u(x).\eeq
It is well known that
(\ref{1}) generates a Lipschitz continuous semigroup
of Liu-admissible weak solutions \cite{BiB, B1, Bbook, Btut, BC95, BCP, BLY, Daf, HR}, on a
domain  $\D\subset \L^1(\R;\,\R^n)$ of functions with suitably small total variation.
In the papers \cite{BC2, Lew1, Lew2, LewT} a semigroup was constructed on a domain of functions with large
but finite total variation. In view of the recent uniqueness theorems proved in \cite{BDL, BGu},  these results show that, for a wide class of systems,  the Cauchy problem (\ref{1})-(\ref{2}) has a unique solution
which depends continuously on the initial data, as long as the
total variation remains bounded.

A major open problem is to understand whether  the semigroup generated
by (\ref{1}) can be extended to a domain of $\L^\infty$ functions, 
possibly with unbounded total variation.
At present, this is known in the
scalar case \cite{Crandall, Kru}, and for some special systems of Temple class \cite{BG2}.  In addition, unique solutions with $\L^\infty$ initial data  were constructed
in \cite{BGS} as limit of vanishing viscosity approximations, 
for a class of $2\times 2$ systems of conservation laws in triangular form, 
not necessarily strictly hyperbolic.

Even for $2\times 2$ strictly hyperbolic, genuinely nonlinear systems,
where existence of $\L^\infty$ solutions can be proved by a compensated compactness
argument \cite{DiP}, it is not known whether the vanishing viscosity limit yields unique solutions. 
In particular, the uniqueness of the weak solutions
constructed by Glimm and Lax in \cite{GL} has remained an elusive open problem
for over 50 years.  

Aim of this paper is to provide some results in this direction.   Our approach is motivated by the theory of parabolic equations and fractional domains for 
analytic semigroups \cite{Henry, Lu}, which we briefly review.  Consider  a bounded open set $\Omega\subset\R^n$  with smooth boundary.
Call $u(t,\cdot) = e^{t\Delta} \ov u$ the solution to Cauchy problem for the linear heat equation
\bel{heateq} \left\{ \bega{rll} u_t - \Delta u &=~0\qquad &\hbox{on}~\Omega,\\[2mm]
u&=~0\qquad  &\hbox{on}~\partial \Omega,\enda\right.\qquad\qquad u(0,x)\,=\,\ov u(x).
\eeq
For any  initial data $\ov u\in \L^2(\Omega)$, it is well known that this solution
satisfies the decay property
\bel{regp} \bigl\|\Delta u(t,\cdot)\bigr\|_{\L^2}~\leq~{C\over t} \,  \,\|\ov u\|_{\L^2}\,.
\eeq
Moreover, if the initial data lies in the domain of a fractional power
of the Laplace operator, say $\bar u\in X^\alpha\doteq Dom (-\Delta)^\alpha$ for some $0<\alpha<1$,
then a faster decay rate is achieved:
\bel{regpa} \bigl\|\Delta u(t,\cdot)\bigr\|_{\L^2(\Omega)}~\leq~{C\over t^{1-\alpha}}\, \|\ov u\|_{X^\alpha}\,.
\eeq
Relying on (\ref{regpa}) one can show that
 the semilinear parabolic Cauchy problem
\bel{pareq} \left\{ \bega{rll} u_t - \Delta u &=~F(x,u, \nabla u)\qquad &\hbox{on}~\Omega,\\[2mm]
u&=~0\qquad  &\hbox{on}~\partial \Omega,\enda\right.\qquad\qquad u(0,x)\,=\,\ov u(x),
\eeq
with $F$ Lipschitz continuous,
is well posed for initial data $\ov u\in X^\alpha$, for a suitable $0<\alpha<1$ depending on the dimension $n$ of the space. For all details we refer to \cite{Henry}.

We shall pursue a similar approach in connection with a $2\times 2$
 genuinely nonlinear,
strictly hyperbolic system of conservation laws, as considered in the classical 
paper \cite{GL}.
We recall that, for any bounded measurable initial data $\ov u:\R\mapsto\R^2$, 
if the norm 
\bel{Linfty}\|\ov u\|_{\L^\infty}~\leq~\ve_1\eeq is sufficiently small, then the Cauchy
problem (\ref{1})-(\ref{2}) has a weak solution, global in time. Its $\L^\infty$ norm and its total variation satisfy  bounds of the form
\bel{Linorm} \bigl\|u(t,\cdot)\bigr\|_{\L^\infty}~\leq~M_0\, \sqrt{\|\bar u\|_{\L^\infty}}\,,\eeq
\bel{TVdec}
\TV\bigl\{ u(t,\cdot)\,;~[a,b]\bigr\}~\leq~C_0\cdot\left(\sqrt{\|\bar u\|_{\L^\infty}} + {b-a\over t}\right)\,,\eeq
for some constants $M_0, C_0$, any interval $[a,b]$ and every time $t\in \,]0, 1]$.    See also \cite{BCM}
for a simpler proof, based on front tracking approximations.

In analogy with (\ref{regpa}), we introduce the domain
\bel{TDA}\bega{rl} \Tilde \D_\alpha&\ds\doteq~\Bigg\{ \ov u\in\L^\infty(\R;\,\R^2)\,; ~\hbox{ the Cauchy problem (\ref{1})-(\ref{2}) has a weak solution}
\\[4mm]
&\qquad\qquad\qquad  \ds \hbox{such that}\quad
\sup_{0<t<1}  t^{1-\alpha} \cdot \TV\bigl\{ u(t,\cdot)\bigr\}\,<\,+\infty
\bigg\} \,.\enda\eeq
We observe that $\bar u\in \Tilde\D_\alpha$ implies
\bel{dista}
\bigl\| u(\tau)-\ov u\|_{\L^1}~=~\O(1)\cdot \tau^\alpha\qquad\qquad \tau\in [0,1].\eeq
Indeed, by Theorem~4.3.1 in \cite{Daf}, one has 
\bel{lis}\limsup_{h\to 0+}~{\bigl \| u(t+h) - u(t)\bigr\|_{\L^1}\over h}~=~\O(1)\cdot \TV\bigl\{ u(t)\bigr\}
\,.\eeq
For $\bar u\in \Tilde\D_\alpha$, the estimate (\ref{dista}) is an easy consequence of (\ref{lis}).
Here and in the sequel, the Landau notation $\O(1)$ denotes a uniformly bounded quantity.

Our main result shows that, for every initial data $\ov u\in \Tilde \D_\alpha$ with $\|\bar u\|_{\L^\infty}$ sufficiently small,
the entropy weak solution which satisfies (\ref{dista})  is unique.  Moreover, 
such solutions form a H\"older continuous semigroup w.r.t.~the $\L^1$ distance on the initial data.
By restricting the domain of the semigroup to functions with even smaller $\L^\infty$ norm, the H\"older exponent of 
the semigroup can be rendered arbitrarily close to 1.

\begin{figure}[ht]
  \centerline{\hbox{\includegraphics[width=10cm]{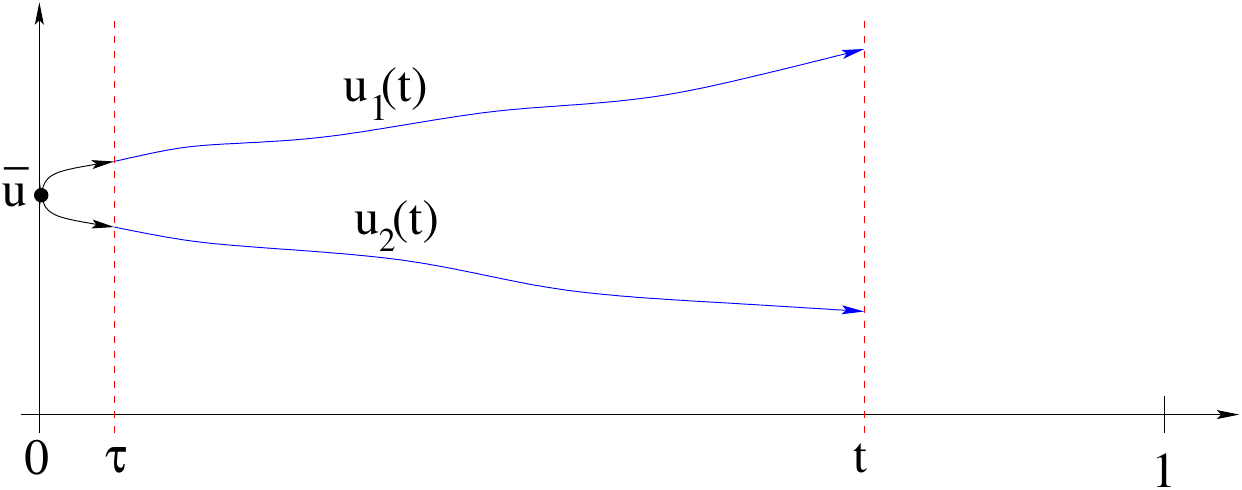}}}
  \caption{\small Sketch of the estimates (\ref{Lip3})--(\ref{es2}). }\label{f:ag49}
\end{figure}

The main idea is illustrated in Fig.~\ref{f:ag49}.  
For any $\tau>0$, consider the positively invariant domain
\bel{Dtau} \bega{l} \D^{(\tau)}~=~\Big\{ u(t,\cdot)\,;~~t\geq \tau,~~u~\hbox{is an entropy weak solution of
(\ref{1})-(\ref{2}) satisfying (\ref{TVdec}),}\\[3mm]
\qquad\qquad\qquad\qquad \hbox{with}~ \bar u\in \L^1(\R;\,\R^2),~ \|\bar u\|_{\L^\infty}\leq \ve_1\Big\}.\enda\eeq
Here $\ve_1>0$ is a fixed constant, small enough so  that the Glimm-Lax decay estimates 
(\ref{TVdec}) on the total variation are valid.  Observing that these functions have uniformly bounded variation on bounded intervals, our first goal is to show that (\ref{1}) generates a Lipschitz semigroup $S$ on the domain
$\D^{(\tau)}$. More precisely
\bel{Lip3}
\bigl\| S_t\bar u-S_t\bar v\bigr\|_{\L^1}~\leq~L{(\tau)}\cdot \|\bar u - \bar v\|_{\L^1}\,,\qquad\qquad 
\forall \bar u,\bar v\in \D^{(\tau)}, ~~t\in [0,1 ],\eeq
for some Lipschitz constant $L(\tau)$.   Of course, since the total variation of functions
$\bar u\in \D^{(\tau)}$ becomes arbitrarily large as $\tau\to 0$, we expect that 
$\lim_{\tau\to 0} L(\tau)=+\infty$.  

Next, consider two solutions $u_1,u_2$ of the 
hyperbolic system
(\ref{1}) with the same initial data $\bar u\in \Tilde D_\alpha$.  For any $0<\tau<t\leq 1$,
by (\ref{Lip3}) and (\ref{dista}),  their distance
satisfies
\bel{es2}\bega{l} \bigl\|u_1(t)-u_2(t)\|_{\L^1}~\leq~L(\tau)\cdot \bigl\|u_1(\tau)-u_2(\tau)\|_{\L^1}\\[2mm]
\qquad
\leq~ L(\tau) \cdot \Big(  \bigl\|u_1(\tau)-\bar  u\|_{\L^1}
+  \bigl\|u_2(\tau)-\bar  u\|_{\L^1}\Big)~=~L(\tau)\cdot \O(1)\cdot\tau^\alpha.\enda\eeq
Keeping $t>0$ fixed and letting $\tau\to 0$, the condition
\bel{aslip}\lim_{\tau\to 0} ~\tau^\alpha \, L(\tau)~=~0\eeq
implies $u_1(t)=u_2(t)$ for all $t\in [0,1]$, proving uniqueness.

To carry out this program, two main ingredients are needed:
\begi
\item[{\bf (1)}]  Estimate the Lipschitz constant $L(\tau)$.
In particular, give an asymptotic estimate on how it grows as $\tau\to 0$.
\item[{\bf (2)}] 
Identify a class of initial data $\ov u$ for which (\ref{dista}) holds.   
As already remarked, this will be the case if $\ov u\in \Tilde \D_\alpha$.
\endi

\begin{figure}[ht]
  \centerline{\hbox{\includegraphics[width=8cm]{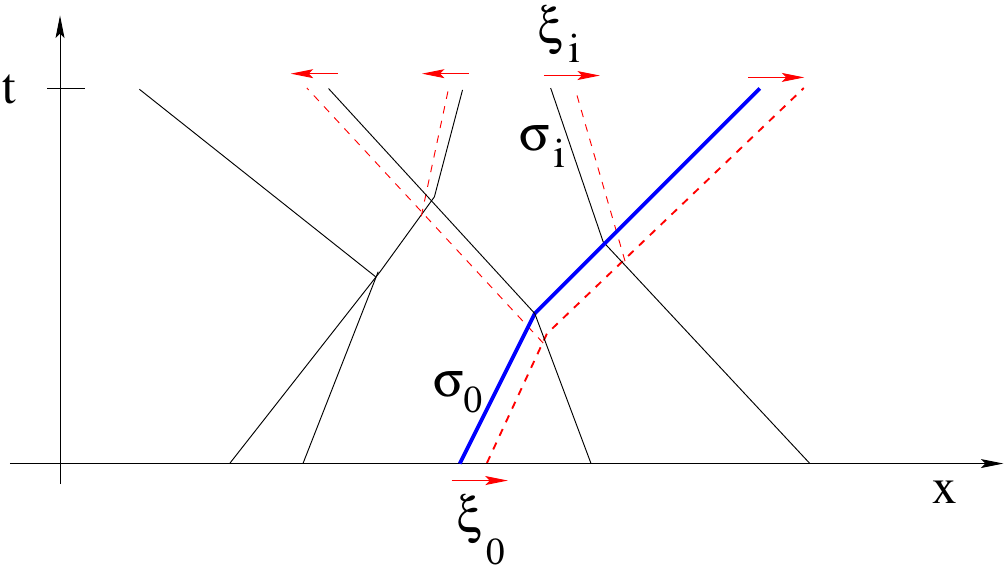}}}
  \caption{\small The Lipschitz continuous dependence of front tracking approximations can be estimated 
  by computing  the strengths $\sigma_i$ and shifts $\xi_i$ of all wave fronts. }\label{f:ag68}
\end{figure}

For part  {\bf (1)}, to estimate the Lipschitz constant of the semigroup we use the technique introduced in 
\cite{BC95}, based on the analysis of strengths and shifts of wave fronts in a front tracking approximation. 
As shown in Fig.~\ref{f:ag68}, consider one particular front at time $t=\tau$, say with size $\sigma_0$.
If the location of this front is initially shifted at rate $\xi_0$, at later time $t>\tau$  a number of other fronts
will be shifted as well.  To fix ideas, call $\sigma_i,\xi_i$, $i=1,\ldots,N$ the sizes and the shifts 
of these various fronts.  
If we can find a constant $L$, uniformly valid for all front tracking solutions, such that
\bel{Lipes}
\sum_i \bigl|\sigma_i \xi_i\bigr|~\leq~L\cdot \bigl|\sigma_0 \xi_0\bigr|,\eeq
this will yield a Lipschitz constant for the semigroup.
We recall that the estimate  (\ref{Lipes}) was first proved in \cite{BC95}, within a domain
of initial data having small total variation.  In the special case of 
genuinely nonlinear Temple class systems \cite{T}, the paper \cite{BG2} showed that the same estimate is valid
on a much larger domain of $\L^\infty$ initial data, which may well have unbounded variation.
A simpler proof of this result for Temple class systems will be given in Section~\ref{sec:3}.\\
A key idea used in our analysis is that any $2\times 2$ system can be locally approximated by a Temple class system.
For a front tracking solution whose $\L^\infty$ norm remains suitably small, the strengths and shifts of wave fronts are thus very close to the corresponding ones computed in connection with  the Temple system. 
This comparison argument will provide the crucial estimate for the Lipschitz constant $L(\tau)$ of our semigroup.

Concerning {\bf (2)} we recall that,
for a scalar conservation law with convex flux, conditions that guarantee a fast decay 
of the total variation have been provided
in the recent paper \cite{ABMT}.  In Section~\ref{sec:8} we show that a strengthened version of these
conditions imply $\ov u\in \Tilde \D_\alpha$, in the case of $2\times 2$ systems.   
More precisely, we consider the set of functions  
$\ov u\in \L^\infty(\R;\,\R^2)$ with the property:
\begi
\item[{\bf ($\Tilde {\bf P}_\alpha$)}]  {\it  
For any  $\lambda\in \,]0,1]$ and every bounded interval $I$, there exists an open set 
$V^\lambda\,=\,\cup_{i=1}^N ]a_i, b_i[\,$ such that the following holds.
\bel{mop2}
\meas(V^\lambda) ~\leq~\Tilde C\, \lambda^\alpha,
\eeq
\bel{tvp2}
\TV\bigl\{ \ov u\,; ~I \setminus V^\lambda\bigr\}~\leq~\Tilde C \,\lambda^{\alpha-1},\eeq
\bel{Nb}
N~\leq~\Tilde C\, \lambda^{\alpha-1},\eeq
for some constant $\Tilde C$ independent of $\lambda$.}
\endi
Roughly speaking, this means that the total variation of $\ov u$ can be infinite,
but most of this variation is concentrated on finitely many open intervals, 
with small total length.

Our main results can be now be stated as

\begin{theorem} \label{t:11} {\bf (Uniqueness for a class of $\L^\infty$ solutions with fast decay).}
Let (\ref{1}) be a $2\times 2$ strictly hyperbolic, genuinely nonlinear system of conservation laws with flux $f\in \C^3$.
Given any constants $\alpha\in \,]0,1[$ and $C>0$,  there exists $\ve_1>0$ such that the following holds.
Assume  $\|\bar u\|_{\L^\infty}<\ve_1$  and let $t\mapsto u(t,\cdot)$ be a weak solution  of (\ref{1})-(\ref{2})
which satisfies the Lax admissibility conditions at every point of approximate jump, together with 
the fast decay estimate
\bel{decal}\sup_{0<t<1} t^{1-\alpha} \cdot \TV\bigl\{  u(t)\bigr\}\, <\, C.\eeq
Then such solution is unique.
Indeed, the family of all Lax-admissible weak solutions that satisfy (\ref{Linfty}) and (\ref{decal})
yields a H\"older continuous semigroup:
\bel{hold} \bigl\| u_1(t, \cdot) - u_2(t,\cdot)\bigr\|_{\L^1} ~\leq~\O(1)\cdot \|\bar u_1 - \bar u_2\|^{\beta}_{\L^1}
\qquad\forall t\in [0,1].\eeq 
Here the exponent $\beta<1$  can be taken arbitrarily close to $1$, by choosing $\ve_1>0$ small enough.   
\end{theorem}

\begin{theorem}\label{t:12}{\bf (Solutions with fast decay rate).} Let (\ref{1}) be a $2\times 2$ strictly hyperbolic, genuinely nonlinear   system of conservation laws, and let $\alpha\in \,]0,1[\,$ be given. Then there exists
$\ve_1>0$ small enough such that the following holds.
Let $\ov u:\R\mapsto \R^2$ be a bounded measurable function
with $\|\ov u\|_{ \L^\infty} <\ve_1$.  If $\ov u$ satisfies the assumptions 
{\bf ($\bf \Tilde P_\alpha$)}, then  the Cauchy problem 
(\ref{1})-(\ref{2}) has a solution $t\mapsto u(t,\cdot)$ with the following decay rate on every bounded interval $[a,b]$:
\bel{decrat}
\TV\bigl\{ u(t,\cdot);\, [a,b]\bigr\}~\leq~C\, t^{\alpha-1}\qquad\qquad\forall t\in \,]0,1].\eeq
\end{theorem}

We remark that  the constants $\ve_1$ and $C$ in (\ref{decrat}) depend on $f$ and on the constants 
$\alpha, \Tilde C$
arising in the statement of  {\bf ($\Tilde {\bf P}_\alpha$)} on the interval $I=\bigl[a - t, b +  t\bigr]$.

The remainder of the paper is organized as follows. In Section~\ref{sec:2} we review 
some basic estimates
on the strengths and shifts of wave fronts across an interaction.  
Section~\ref{sec:3} briefly recalls
the construction of piecewise constant front tracking solutions introduced in \cite{BC95}.
We remark that the original front tracking algorithm in \cite{B92} applied to general $n\times n$ 
systems, but required the introduction of ``nonphysical wave-fronts" to prevent the total number of fronts from blowing up in finite time.  On the other hand, for $2\times 2$ systems, the algorithm introduced in \cite{BC95} uses a
smooth interpolation between shock and rarefaction curves and completely avoids the
introduction of non-physical fronts.  This is a crucial property, because it 
guarantees that front tracking approximations depend continuously on the initial data.
The $\L^1$ distance between two front tracking solutions at time $t>0$ can now be estimated by a homotopy method.
Given two piecewise constant initial data $\bar u^0$, $\bar u^1$,  by shifting the positions of the initial jumps 
one constructs a 1-parameter family of initial data $\bar u^\theta$, $\theta\in [0,1]$.   At a later time $t$, 
the distance $\bigl\|u^0(t,\cdot) -u^1(t,\cdot)\bigr\|_{\L^1}$ is bounded by the length of the path
$\theta\mapsto u^\theta(t,\cdot)$.
In turn (see Fig.~\ref{f:ag68}), this length can be bounded by estimating  the strengths and shifts of all wave-fronts in $u^\theta$.

In the present paper, by covering the $t$-$x$ plane with trapezoidal domains where the total
variation remains uniformly small, from the results in \cite{BC95} we obtain
a preliminary estimate on the total amount of interaction and on the shifts
of all wave-fronts, in a front tracking approximation.  This yields a first, rough estimate of the Lipschitz dependence on the initial data $\bar u\in\D^{(\tau)}$.

Section~\ref{sec:4} is concerned with 
a $2\times 2$ genuinely nonlinear Temple class system, with $\L^\infty$ initial data.
We show that the entropy admissible solutions yield a continuous semigroup, with uniformly Lipschitz dependence on the initial data, in the $\L^1$ distance. This result was originally proved in \cite{BG2}.
We provide here a more direct proof, which yields a sharper Lipschitz constant.

In Section~\ref{sec:5} we show that any $2\times 2$ genuinely nonlinear system can be locally
approximated by a Temple system. Thanks to this approximation, in Section~\ref{sec:6}
we provide a new, sharper bound on the strengths and shifts
for a front tracking solution of the original system (\ref{1}).
To explain the heart of the matter we observe that, for a solution with small $\L^\infty$ norm,
genuine nonlinearity produces a huge cancellation of positive and negative waves. This is indeed
exploited in the famous Glimm-Lax paper \cite{GL}.  One expects that a similar cancellation
should occur also for the shifts in a front tracking solution.  
This can be clearly seen for Temple systems, but is far from obvious in the case of more general 
systems.  Indeed, the bounds in \cite{BC95} do not account for these cancellations.
Our strategy is to approximate an arbitrary $2\times 2$ system with a Temple system, and 
compare the strengths and shifts of wave fronts in the corresponding front tracking solutions.
It turns out that this difference can be rendered arbitrarily small by reducing the
$\L^\infty$ norm of the initial data. 

Having collected all these preliminary estimates,  in Section~\ref{sec:7} we give a proof of Theorem~\ref{t:11}, while 
Section~\ref{sec:8} contains a proof of Theorem~\ref{t:12}.

Additional results on existence and stability of solutions to balance laws with possibly unbounded variation
can be found in  \cite{ABBCN, BS1, HJ, Lew99}. For a comprehensive introduction to hyperbolic conservation laws
we refer to \cite{Bbook, Daf, HR, Liu}.

\section{Estimates on the strengths and shifts of wave fronts} 
\label{sec:2}
\setcounter{equation}{0}
Following the approach introduced in \cite{BC95}, the Lipschitz constant $L(\tau)$ of the 
semigroup generated by (\ref{1}) on the domain $\D^{(\tau)}$ will be estimated by 
bounding strengths and shifts of wave-fronts, in a front tracking approximation. 
In this preliminary section we review, and slightly refine, 
the interaction estimates proved in \cite{BC95}.

Let $(w_1, w_2)$ be 
a set of Riemann coordinates  for the system (\ref{1}). 
In terms of these coordinates, for smooth solutions the system takes a diagonal form:
\bel{diag}\left\{ \bega{rl}w_{1,t}+\lambda_1(w_1, w_2) w_{1,x}&=~0,\\[1mm]
w_{2,t}+\lambda_2(w_1,w_2) w_{1,x}&=~0.\enda\right.
\eeq
Consider an elementary wave (either a rarefaction or a shock) with left and right
states
\bel{lrstate}u^-~=~u(w_1^-, w_2^-),\qquad\qquad u^+~=~u(w_1^+, w_2^+). \eeq
If the wave is of the $i$-th family, $i\in \{1,2\}$,
the size of this jump is defined to be
$$\sigma_i~\doteq~w_i^+-w_i^-.$$
The next lemma collects the basic estimates on the strength and shifts of wave fronts, across an interaction.
\begin{figure}[ht]
  \centerline{\hbox{\includegraphics[width=10cm]{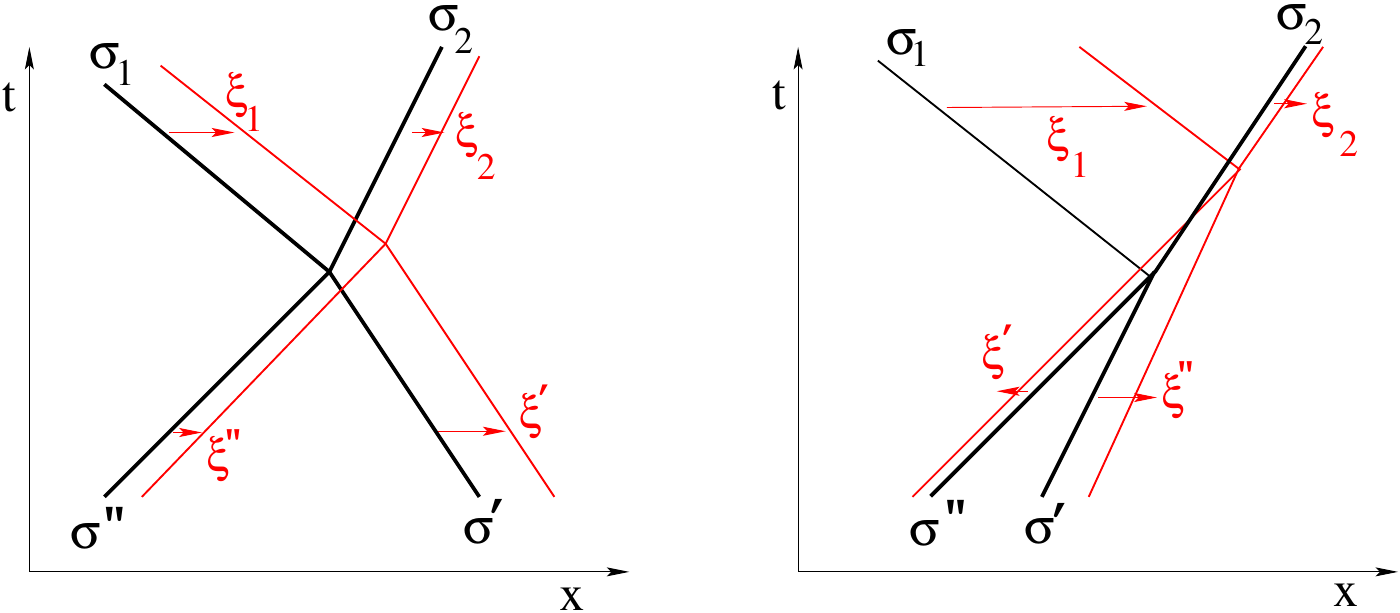}}}
  \caption{\small  The change in strength and shift 
  of two outgoing fronts, depending on the strength, shifts, and relative speeds of the incoming fronts. }
  \label{f:ag51}
\end{figure}
\begin{lemma}\label{l:21}
Consider an interaction of two incoming fronts with sizes $\sigma', \sigma''$, and 
let $\xi', \xi''$ be the corresponding shifts
(see Fig.~\ref{f:ag51}).

Then the sizes of the outgoing fronts $\sigma_1, \sigma_2$
and their shifts $\xi_1,\xi_2$
satisfy the following  estimates.
\v
CASE 1: If the two incoming fronts belong to different families, then
\bel{ies1}
|\sigma_1 -\sigma'| + |\sigma_2-\sigma''|~=~\O(1) \cdot |\sigma'\sigma''| \bigl(|\sigma'|^2+|\sigma''|^2\bigr),\eeq
\bel{ies2}
|\xi_1 -\xi'|~=~\O(1) \cdot |\sigma''|\bigl( |\xi'|+|\xi''|\bigr), \qquad\qquad  |\xi_2 -\xi''|~=~\O(1) \cdot |\sigma'|\bigl( |\xi'|+|\xi''|\bigr).\eeq

CASE 2: If the two incoming fronts belong to the same  family (say, the first one) then
\bel{ies3}
|\sigma_1 -\sigma'-\sigma''|+|\sigma_2| ~=~\O(1) \cdot |\sigma'\sigma''| \bigl(|\sigma'|+|\sigma''|\bigr),\eeq
\bel{ies4}~\left| \xi_1-{\xi'\sigma'+\xi''\sigma''\over \sigma'+\sigma''} \right|~=~\O(1)\cdot \bigl( |\xi'|+|\xi''|\bigr) {|\sigma'\sigma''| \over  {|\sigma'|+|\sigma''|}},
\qquad\qquad 
|\xi_2|~=~\O(1) \cdot {|\xi'|+|\xi''|\over |\sigma'|+|\sigma''|}.\eeq
\end{lemma}
\v

\begin{figure}[ht]
  \centerline{\hbox{\includegraphics[width=12cm]{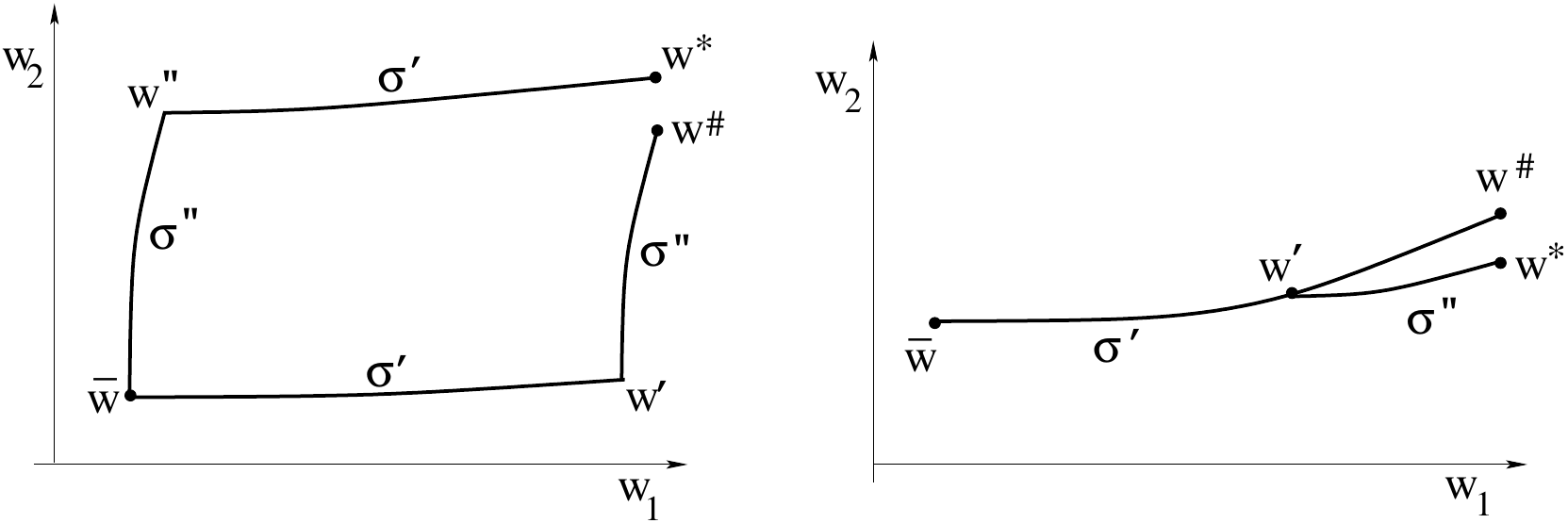}}}
  \caption{\small Estimating the sizes of wave fronts, before and after an interaction.}\label{f:ag56}
\end{figure}

{\bf Proof.} {\bf 1.}  To prove (\ref{ies1}), fix a state $\ov w=(\ov w_1,\ov w_2)$.
Then the first and second shock curves through $\ov w$ are described by
$$s~\mapsto~S_1(s, \ov w)~=~\Big( \ov w_1+s,~\ov w_2 + \vp( s,\ov w )\, s^3\Big),$$
$$s~\mapsto~S_2(s, \ov w)~=~\Big( \ov w_1+\psi(s,\ov w)\,s^3,~\ov w_2 +s\Big),$$
for some smooth functions $\vp,\psi$.
Consider two incoming shocks of different families,
of size $\sigma'$ and $\sigma''$, respectively.   Calling 
$\ov w, w'', w^*$ the left, middle, and right state before interaction, as shown in Fig.~\ref{f:ag56}, left,  we have
$$w''~=~S_2(\sigma'', \ov w)~=~\Big( \ov w_1+\psi(\sigma'',\ov w)\cdot (\sigma'')^3,~\ov w_2 +\sigma''\Big),$$
$$w^*~=~S_1(\sigma', w'')~=~\Big( w''_1+\sigma',~w_2'' + \vp( \sigma' ,w'' )\cdot (\sigma')^3\Big).$$
Next, consider the points 
$$w'~=~S_1(\sigma', \ov w)~=~\Big( \ov w_1+\sigma',~\ov w_2 + \vp( \sigma' ,\ov w )\cdot (\sigma')^3\Big).$$
$$w^\sharp~=~S_2(\sigma'',  w')~=~\Big( w'_1+\psi(\sigma'', w')\cdot (\sigma'')^3,~ w'_2 +\sigma''\Big),$$
By the implicit function theorem we have
\bel{ies 6}\bega{l}  |\sigma_1-\sigma'| + |\sigma_2 - \sigma''| 
~=~\O(1) \cdot |w^\sharp- w^*|\\[4mm]
\qquad =~\O(1)\cdot \Big| \psi(\sigma'', \ov w) -\psi(\sigma'', w')\Big| | \sigma''|^3 +
\O(1) \cdot \Big| \psi(\sigma', \ov w) -\psi(\sigma'', w'')\Big| | \sigma'|^3\\[4mm]
\qquad =~\O(1) \cdot |\sigma'|\, \bigl|\sigma''\bigr|^3 + \O(1)\cdot  |\sigma''|\, \bigl|\sigma'\bigr|^3.
\enda
\eeq
This yields (\ref{ies3}) in the case where both interacting fronts are shocks.

The case where one of the incoming fronts is a rarefaction is treated in a similar way.
If both incoming waves are rarefactions, we trivially have $\sigma_1=\sigma'$, $\sigma_2=\sigma''$.
\v
{\bf 2.} Next, consider the interaction of two shocks of the same family,  say  of the first one.
Let $\ov w, w', w^*$ be the left, middle, and right states before the interaction.
We thus have
$$w'~=~S_1(\ov w, \sigma')~=~\Big( \ov w_1 + \sigma',~\ov w_2 + \vp(\ov w, \sigma') (\sigma')^3\Big).$$
$$w^*~=~S_1(w', \sigma'')~=~\Big(  w_1' + \sigma'',~w'_2 + \vp( w', \sigma'') (\sigma'')^3\Big).$$
Consider the point
$$w^\sharp~=~S_1(\ov w, \sigma'+\sigma'')~=~\Big( \ov w_1 + \sigma'+\sigma'',~\ov w_2 + \vp(\ov w, \sigma'+\sigma'') (\sigma'+\sigma'')^3\Big).$$
By the implicit function theorem we have
\bel{ies7} \bega{l}
|\sigma_1 -\sigma'-\sigma''| + |\sigma_2|~=~\O(1)\cdot |w^\sharp-w^*|\\[3mm]
\qquad =~\O(1)\cdot \Big|  \vp( \ov w, \sigma'+\sigma'') (\sigma'+\sigma'')^3- \vp(\ov w, \sigma') (\sigma')^3 - \vp( w', \sigma'') (\sigma'')^3\Big|\\[3mm]
\qquad =~\O(1) \cdot \bigl|\sigma'\bigr|^2 |\sigma''| + \O(1)\cdot |\sigma'|\bigl|\sigma''\bigr|^2.
\enda\eeq
This proves (\ref{ies3}).
\v
{\bf 3.} In the remainder of the proof we study the shifts of the outgoing fronts.
Calling $\lambda'< \lambda''$ the speeds of the incoming fronts and $\lambda_1, \lambda_2$ 
the speeds of the outgoing fronts of the first and second family, the interaction time changes by
$$\Delta t~=~{\xi'-\xi''\over \lambda''-\lambda'}\,.$$
In turn, the shift in the outgoing front of the first family is computed by
\bel{shi1}\xi_1~=~\xi' -(\lambda_1-\lambda') \Delta t~=~{\xi'(\lambda''-\lambda_1) -\xi''(\lambda'-\lambda_1)\over \lambda''-\lambda'}\,.\eeq
Similarly,  the shift in the outgoing front of the second family is computed by
\bel{shi2}\xi_2~=~\xi'' -(\lambda_2-\lambda'') \Delta t~=~{\xi'(\lambda''-\lambda_2) -\xi''(\lambda'-\lambda_2)\over \lambda''-\lambda'}\,.\eeq
\v
{\bf 4.} Consider two interacting fronts of distinct families.
In this case, the change in speeds across the interaction is estimated by
\bel{shi33}|\lambda_1-\lambda'|~=~\O(1)\cdot |\sigma''|, \qquad\qquad |\lambda_2-\lambda''|~=~\O(1)\cdot |\sigma'|.\eeq 
Since the difference in speeds $\lambda''-\lambda'$ is bounded away from zero, 
by (\ref{shi1})-(\ref{shi2}) this yields (\ref{ies2}).
\v
{\bf 5.}  Finally, consider two interacting fronts of the same family, say, the first one.
In this case, by genuine nonlinearity, the difference in speed has the lower bound
\bel{ies6}|\lambda'-\lambda''|~\geq~c_0 \bigl(|\sigma'|+|\sigma''|\bigr).\eeq
This bound is clear if both incoming fronts are shocks. If one is a rarefaction, we
can assume that its size is smaller that 1/2 times the size of the shock,
so that (\ref{ies6}) still holds. By (\ref{shi2}) this implies
the second estimate in (\ref{ies4}).

On the other hand, the speed $\lambda_1$ of the outgoing front of the first family
satisfies
\bel{ospe}
\left| \lambda_1-{\sigma'\lambda' + \sigma''\lambda''\over \sigma'+\sigma''}\right|
~=~\O(1)\cdot |\sigma'\sigma''|.\eeq
Inserting (\ref{ospe}) in (\ref{shi1}) one obtains
\bel{ies71}\bega{rl} \xi_1&\ds=~{1\over \lambda''-\lambda'}\cdot \left\{
\xi' \left({ (\lambda''-\lambda')\sigma' \over  \sigma'+\sigma''}  \right) 
+\xi'' \left({(\lambda '' -\lambda')\sigma''\over \sigma'+\sigma''}\right)+
\O(1)\cdot \bigl( |\xi'|+|\xi''|\bigr) |\sigma'\sigma''|\right\}\\[4mm]
&\ds =~{\xi'\sigma'+\xi''\sigma''\over \sigma'+\sigma''} +\O(1)\cdot \bigl( |\xi'|+|\xi''|\bigr) {|\sigma'\sigma''| \over |\sigma'|+|\sigma''|}.\enda\eeq
This completes the proof. \endproof

\begin{cor} In an interaction of fronts of different families, the increase in the total amount of
wave strengths is bounded by
\bel{ies 60}|\sigma_1|+|\sigma_2| - |\sigma'|-|\sigma''|~\leq~\O(1)\cdot \, |\sigma'\sigma''|
 \bigl(|\sigma'|^2+|\sigma''|^2
\bigr).\eeq 
If the incoming fronts are of the same family, one has 
\bel{ies 61}|\sigma_1|+|\sigma_2| - |\sigma'|-|\sigma''|~\leq~\O(1)\cdot \, |\sigma'\sigma''| \bigl(|\sigma'|+|\sigma''|
\bigr).\eeq 
In both of the above cases, the shifts satisfy 
\bel{ies72}
|\xi_1 \sigma_1| + |\xi_2\sigma_2| - \bigl( |\xi'\sigma'|+|\xi''\sigma''|\bigr)~\leq~
\O(1)\cdot \bigl(|\xi'|+|\xi''|\bigr)\, |\sigma'\sigma''|.\eeq
\end{cor}

\begin{remark} {\rm If the shifts of the incoming fronts are $\xi'=\xi''=1$, then one 
has $\xi_1=\xi_2=1$ as well.   In this special case, the bounds (\ref{ies72}) reduce to the standard 
Glimm interaction estimates on the strengths of outgoing fronts \cite{Glimm}.
}
\end{remark}
%

\section{Front tracking approximations}
\label{sec:3}
\setcounter{equation}{0}
We review here the construction of front tracking approximations introduced in \cite{BC95}.
Working with a coordinate system of Riemann invariants, consider a point with coordinates 
$w=(w_1, w_2)$.   The rarefaction curves through this point are parameterized as
\bel{rar}s~\mapsto~ R_1(w,s)\,=\,(w_1+s, w_2), \qquad\qquad s~\mapsto~ R_2(w,s)\,=\,(w_1, w_2+s).
\eeq
On the other hand, for $s<0$ the shock curves can be parameterized as
\bel{shock} 
s~\mapsto~ S_1(w,s)\,=\,\bigl(w_1+s, w_2+\phi_2(w,s)\bigr), \quad
\qquad s~\mapsto~ S_2(w,s)\,=\,\bigl(w_1+\phi_1(w,s), w_2+s\bigr),\eeq
where $\phi_1,\phi_2$ are smooth functions with $\phi_i(w,s) = \O(1)\cdot s^3$.

\begin{figure}[ht]
  \centerline{\hbox{\includegraphics[width=7cm]{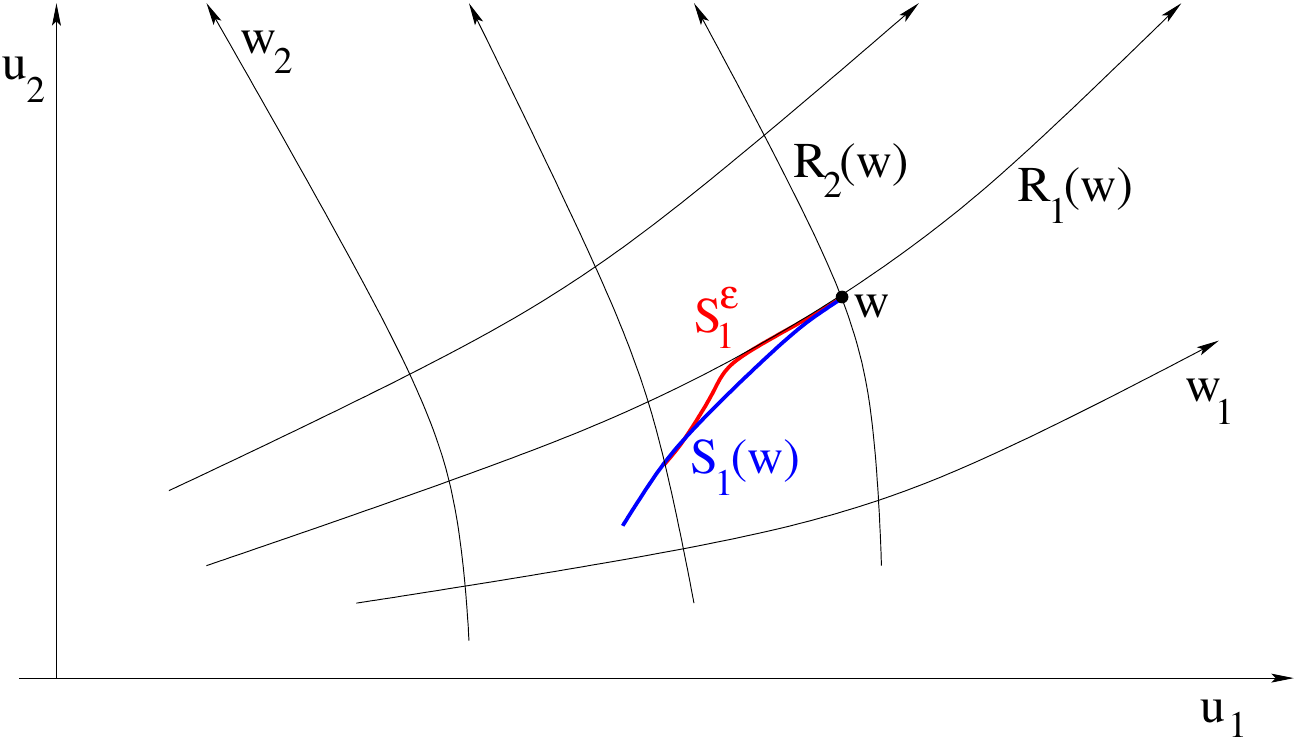}}\qquad \hbox{\includegraphics[width=8cm]{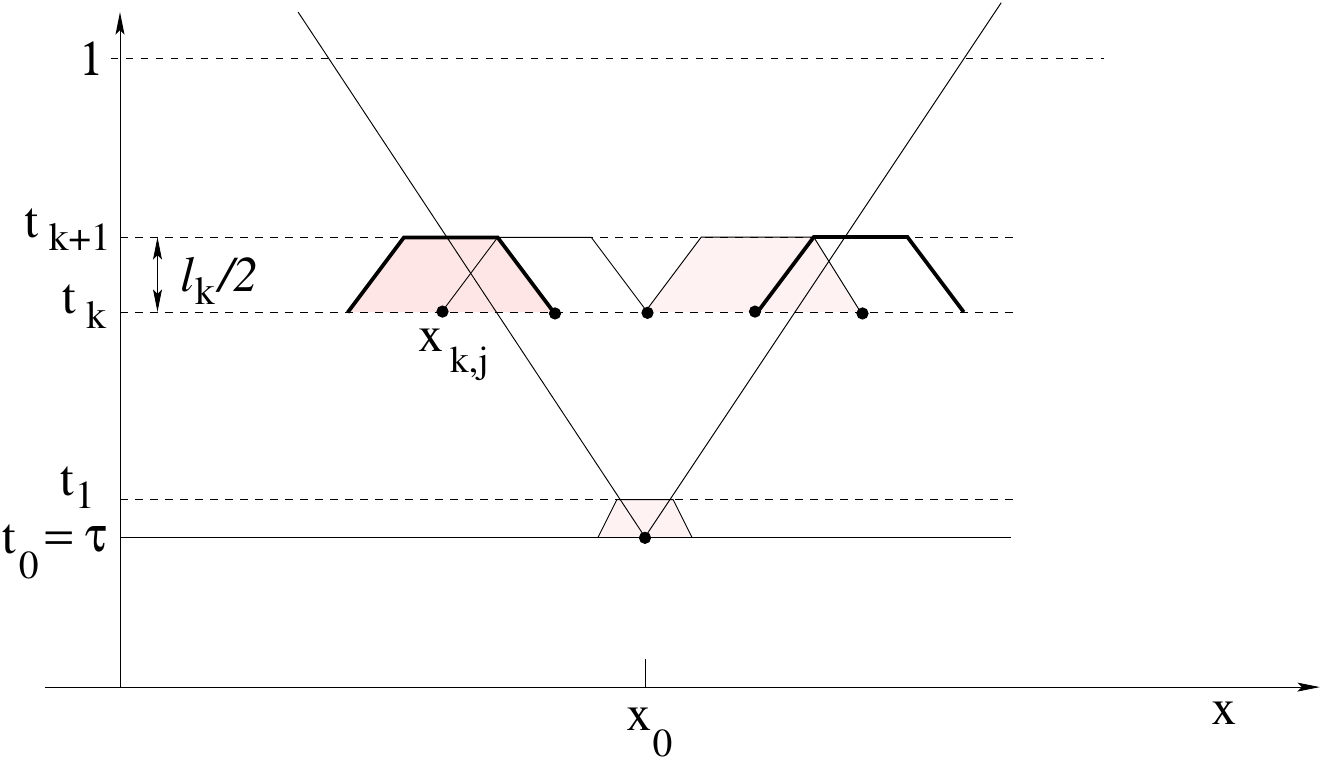}}}
  \caption{\small Left: a system of Riemann coordinates $(w_1,w_2)$ for the system (\ref{1}). Here
  $s\mapsto S_1^\ve(w,s)$ is the interpolated shock curve defined at (\ref{Seps}).  Right:
  covering the triangular region $\Delta$ with trapezoids $\Gamma_{kj}$ where the total variation
is suitably small, by (\ref{txk})-(\ref{tgk}) one obtains a preliminary upper bound on the 
 shifts in a front tracking solution.}\label{f:ag86}
\end{figure}

We now fix $\ve>0$ and define the interpolated shock curves $S_1^\ve, S_2^\ve$ so that
(see Fig.~\ref{f:ag86}, left)
$$S_i^\ve(w,s)~=~\left\{\bega{rl} R_i(w,s)\quad &\hbox{if}\quad -\sqrt\ve<s<0,\\[1mm]
S_i(w,s)\quad &\hbox{if}\quad s<-2\sqrt\ve
.\enda
\right.
$$
This is achieved by constructing a nondecreasing smooth function
$\vp:\R\mapsto [0,1]$ such that 
$$\vp(t)~=~\left\{\bega{rl}0\quad &\hbox{if}\quad t<-2,  \\[1mm]
1\quad &\hbox{if}\quad  t\geq -1,\enda\right.$$
and defining
\bel{Seps}S^\ve_i(w,s)~=~\vp(s/\sqrt \ve) R_i(w, s) + \bigl(1-\vp(s/\sqrt \ve)\bigr) S_i(w,s).\eeq
In turn, we can now define the mixed curves 
\bel{mix}\Psi^\ve_i(w,s)~=~\left\{\bega{rl}R_i(w,s)\quad &\hbox{if}\quad s\geq 0,  \\[1mm]
S_i^\ve(w,s)\quad &\hbox{if}\quad  s<0.\enda\right.\eeq

Given a left and a right state $w^-=(w_1^-, w_2^-)$, $w^+=(w_1^+, w_2^+)$,
sufficiently close to each other, using the fact that shock and rarefaction curves have a second 
order tangency, one proves that there exists $s_1,s_2$ and a unique intermediate state $w^*=(w_1^*, w_2^*)$ such that
$$w^*~=~\Psi_1^\ve(w^-, s_1), \qquad\qquad w^+~=~\Psi_2^\ve(w^*, s_2).$$
A piecewise constant front tracking 
solution of the Riemann problem with left and right states $w^-, w^+$ 
is then constructed by approximating every rarefaction wave with a centered rarefaction fan.
More precisely, a 1-rarefaction wave with left and right states 
$$w^-=(w_1^-, w_2^-), \qquad w^*=R_1(w^-,s) = (w_1^-+s, w_2^-)$$ is approximated
with a rarefaction fan whose intermediate states lie in the set
$$\Big\{ (j\ve, w_2^-)\,;~~w_1^- \,\leq\, j\ve\,\leq\, w_1+s,\qquad j\in {\mathbb Z}\Big\}.$$
A 2-rarefaction is approximated by a centered rarefaction fan in an entirely similar way.
For all details we refer to \cite{BC95}.

Next, let $u=u(t,x)$ be a front tracking solution of (\ref{1})-(\ref{2}), where the initial datum $\bar u$ satisfies the
bound (\ref{Linfty}) with $\ve_1>0$ sufficiently small.   As proved in \cite{BCM}, this implies the decay 
(\ref{TVdec}) of the total variation. Without loss of generality, we always assume that all wave fronts travel with speed $\bigl|\dot x(t)\bigr|<1$.

Assume that at some time $0<\tau<1$ one of the fronts, with strength $|\sigma_0|$, 
is shifted at rate $\xi_0$.  We seek an estimate on the sum
\bel{Ups1}V^\xi~\doteq~\sum_{\alpha} \bigl|\sigma_\alpha \xi_\alpha\bigr|\eeq
of all the shifts of the  fronts at time $t=1$.

Toward this goal we recall that,  the analysis in \cite{BC95}, there exist constants $\delta_0>0$ and 
 $L>0$, depending only on the system (\ref{1}), such that the following Lipschitz continuity property holds. 
\begi
\item[{\bf (LP)}] {\it Assume that at time $t$ the total strength of wave fronts in $u(t,\cdot)$
contained in the open interval $]a,b[$ is $\leq \delta_0$.     Then for any $t'>t$, one has}
\bel{shb1}\sum_{a+(t'-t)\leq x_\beta\leq b-(t'-t)}\bigl|\sigma_\beta(t')\,\xi_\beta(t')\bigr|~
\leq~L\cdot 
\sum_{a<x_\alpha<b} \bigl|\sigma_\alpha(t)\,\xi_\alpha(t)\bigr|.\eeq
\endi
In other words, let the shifts $\xi_\alpha(t)$ be assigned arbitrarily at time $t$, to all fronts
in the interval $]a,b[$.  Then, at the later time $t'$, 
the shifts of all the fronts in the interval of dependency $\bigl[a+(t'-t), ~b-(t'-t)\bigr]$ satisfy the bound (\ref{shb1}).

As shown in Fig.~\ref{f:ag86},  right, an estimate on the quantity $V^\xi$ in (\ref{Ups1}) can be obtained by 
covering the triangular domain
$$\Delta~\doteq~\Big\{ (t,x)\,;~~t\in [\tau, 1]\,,~~x_0 -(t-\tau)\,<\,x\,<\, x_0+ (t-\tau)\Big\}$$
with finitely many trapezoids $\Gamma_{kj}$ where the total variation is $\leq \delta_0$, and using (\ref{shb1}) iteratively.   

More precisely (see Fig.~\ref{f:ag86}), consider a sequence of times
$$\tau\,=\, t_0\,<\, t_1\,<~~\cdots~~<\, t_N\,,
$$
lengths $\ell_k$ and points $x_{kj}$ satisfying the following inductive relations.
\bel{txk}t_{k+1} ~= ~t_{k} + {1\over 2} \ell_{k}\,,\qquad\qquad x_{kj}\,=\, x_0 + j\cdot \ell_k\,, \eeq
\bel{ellk}\ell_k~=~{\delta_0\over C} t_k\,.\eeq
Here the constant $C$ is chosen so that  
\bel{TVk}
[
\hbox{Total strength of wave-fronts in $u(t_k,\cdot)$
over any interval of length $2\ell_k$}] ~\leq ~\delta_0.\eeq   Notice that, such a constant
exists, uniformly valid  for all $t_k>\tau$, because of (\ref{TVdec}).  
Inserting (\ref{ellk}) in (\ref{txk}) we obtain
\bel{tgk} t_{k+1} ~=~t_k + {\delta_0\over 2C} t_k\,,\qquad\qquad t_k~=~e^{\gamma k} t_0\,,\eeq
where the constant $\gamma>0$ is determined by the identity
$$e^\gamma~=~1+{\delta_0\over 2C}\,.$$
Observing that every point in the triangle $\Delta$ lies in at most two of the trapezoids 
\bel{trapz}\Gamma_{kj}~\doteq~
\bigl\{ (t,x)\,;~~t\in [t_k, t_{k+1}[\,,~~x_{kj} -\ell_k+(t-t_k) <x< x_{kj} + \ell_k - (t-t_k)
\bigr\},\eeq
by (\ref{shb1}) for every $k$ we obtain
\bel{shb2}\sum_{|x_\beta(t_{k+1}) - x_0|< t_{k+1}-\tau} \bigl|\sigma_\beta(t_{k+1})\,\xi_\beta(t_{k+1})\bigr|~
\leq~2L\cdot 
\sum_{|x_\alpha(t_{k}) - x_0|< t_{k}-\tau}  \bigl|\sigma_\alpha(t_k)\,\xi_\alpha(t_k)\bigr|.\eeq
By induction, at any time $T\in [\tau,1]$ this yields
\bel{shb3}V^\xi(T)~=~\sum_{|x_\alpha(T) - x_0|< T-\tau}  \bigl|\sigma_\alpha(T)\,\xi_\alpha(T)\bigr|
~\leq~(2L)^N \, |\sigma_0\xi_0| .\eeq
Here the integer $N$ can be recovered from 
\bel{shb4}e^{\gamma N}\tau\,\geq 1\,,\qquad\qquad N \in  ~\left[{|\ln\tau|\over\gamma}\,,~
~{|\ln\tau|\over\gamma}+1\right]\,.\eeq

\begin{figure}[ht]
  \centerline{\hbox{\includegraphics[width=11cm]{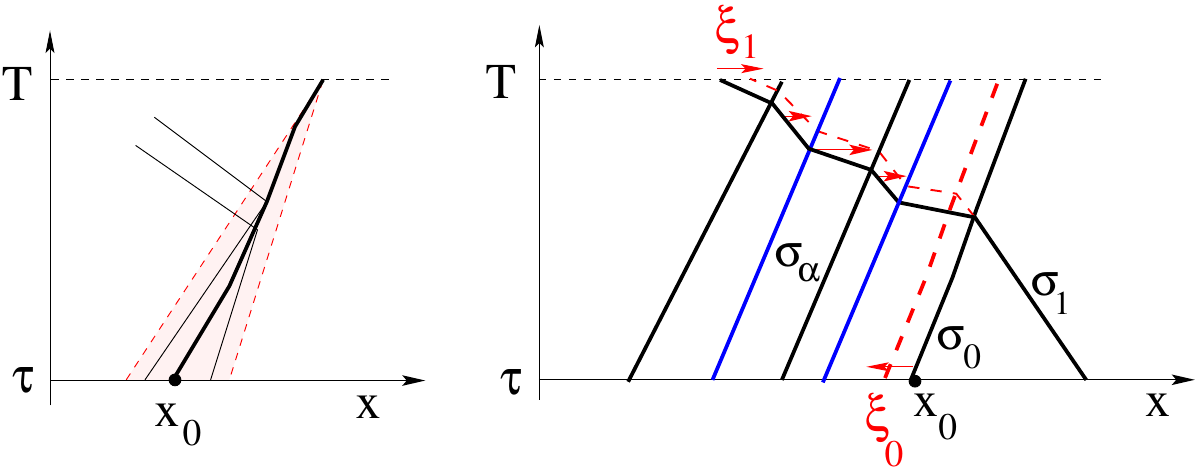}}}
  \caption{\small Left: Since waves of the same family have almost the same speed, the 
  total amount of 2-waves which interact with the 2-wave starting at $x_0$ is small.
  Right: As the 1-front $\sigma_1$ crosses several 2-fronts, its shift $\xi_1$ 
  keeps changing, but a huge cancellation is achieved.}
  \label{f:ag89}
\end{figure}

\begin{remark}\label{r:31} {\rm
The Lipschitz constant $(2L)^N$ in (\ref{shb3})-(\ref{shb4}) is not sharp.
Indeed, it grows very quickly as $\tau\to 0$.   Using the additional information that the $\L^\infty$ norm of the solution is small, a much better estimate will be derived.  The heart of the matter is shown in Fig.~\ref{f:ag89}.
At the initial time $\tau$, the total strength of wave  fronts is  $\O(1)\cdot \tau^{-1}$.   As a consequence,
the total amount of all interactions can be extremely large. However,
the total amount of interaction  {\em restricted to waves of the same family} can be rendered arbitrarily small
by choosing the  initial data small enough in $\L^\infty$.
Indeed (see Fig.~\ref{f:ag89}, left), if $\bigl\| u(t,\cdot)\bigr\|_{\L^\infty}$ remains small, 
then any two wave-fronts of the same family
will have almost the same speed.  Hence they can interact only if their initial locations are very close. 
This will be the content of Lemma~\ref{l:32}.

On the other hand, the amount of interactions among waves of different families can be very large.
In particular (see Fig.~\ref{f:ag89}, right),  if a 1-front crosses a family of 2-fronts, its shift $\xi_1$ keeps changing,
so that its total variation in time is bounded but  possibly very large.  This rough bound is captured by the estimate (\ref{shb3}).   However, since the $\L^\infty$ norm of the solution is small, 
 the sign of the crossing 2-fronts (compressions or rarefactions) keeps alternating. As a result, the shift $\xi_1$
alternatively increases or decreases, and a huge cancellation is achieved.   

Keeping track of this cancellation is not an easy task.   In Section~\ref{sec:4} we refine the arguments in \cite{BaB, BG2} and achieve an explicit, sharp estimate on the strengths and shifts of wave fronts in the special case of a Temple class system.
In the following Sections~\ref{sec:5} and \ref{sec:6}, the strengths and shifts of wave fronts 
in a solution to the original system (\ref{1}) will be estimated by a comparison with a suitable Temple system, relying
on a homotopy method. This will be the content of Lemma~\ref{l:61}. }
\end{remark}

The main result of  this section,  stated in Lemma~\ref{l:32}, 
provides a bound on the shifts produced by interactions of fronts of the same family.
Consider  the total amount of {\it interactions
with shifts}, occurring during the time interval $I_k=[t_k, t_{k+1}[$.  This is defined as
\bel{qshi}\mathcal{Q}_k^\xi~\doteq~\sum_{t\in I_k} \bigl(|\xi_\alpha|+|\xi_\beta|\bigr) |\sigma_\alpha\sigma_\beta|,\eeq
where the sum ranges over all wave-front interactions occurring during the interval $I_k$.
Here $\sigma_\alpha, \sigma_\beta$ and $\xi_\alpha, \xi_\beta$ are the sizes and the shifts of the
two incoming wave fronts, respectively.

Still from the analysis in \cite[Lemma 23]{BC95}, for each $k=0,1,\ldots, N-1$ we have a bound of the form
\bel{Qke}
\mathcal{Q}_k^\xi~\leq ~\delta_1\,V^\xi(t_k)~=~\delta_1\cdot \sum_{|x_\alpha(t_k) - x_0|< t_k -\tau} |\sigma_\alpha(t_k)\xi_\alpha(t_k) \bigr|,\eeq
for some constant $\delta_1$ depending only on  $\delta_0$ in {\bf (LP)} and on the system (\ref{1}).

Summing over $k$ we conclude that the total amount of interaction with shifts satisfies
\bel{Qe}
\mathcal{Q}^\xi~=~\sum_{k=0}^{N-1} \mathcal{Q}_k^\xi~\leq~\delta_1\,\sum_{k=0}^{N-1} (2L)^k |\sigma_0\xi_0|~=~
{(2L)^N-1\over 2L-1} \cdot \delta_1 \,|\sigma_0\xi_0|.\eeq
Notice that, as $\tau\to 0$, by (\ref{shb4}) one has $N\to \infty$ and hence $\Q^\xi\to \infty$ as well.
On the other hand, the next lemma shows that, by choosing $\|\bar u\|_{\L^\infty}<\ve_1$ small enough,
the total amount of interactions with shift, involving fronts of the same family, can be made arbitrarily small.

\begin{lemma}\label{l:32} Given $\tau\in \,]0,1]$, there exists $\ve_1>0$ small enough so that the following holds.   
Let $u=u(t,x)$ be a front tracking solution to (\ref{1})-(\ref{2}), where the initial data
satisfies (\ref{Linfty}). Assume that at time $\tau$ a shift $\xi_0$ is assigned to a front of strength $\sigma_0$.
Then the total amount of interactions with shift, during the interval $t\in [\tau, 1]$ among fronts of the same family,
satisfies the bound
\bel{HatQ}\Q^\xi_{same}~\doteq~\sum_{t\in [\tau,1]\atop {\rm same~ family}} \bigl(|\xi_\alpha|+|\xi_\beta|\bigr) |\sigma_\alpha\sigma_\beta|~\leq~
C_\tau \cdot\sqrt{\ve_1}\,|\sigma_0\xi_0|,\eeq
for some constant $C_\tau$ depending only on $\tau$ and on the system (\ref{1}).
\end{lemma}

{\bf Proof.} {\bf 1.} Let $\Gamma_{kj}$ be any one of the trapezoidal  regions in (\ref{trapz}). 
Restricted to this region, consider the functional
\[
 \Ups_{kj} (t) ~\doteq ~\Big[Q_{same}^\xi(t) + \sqrt\ve_1 \bigl( V^\xi(t) + K Q^\xi(t)\bigr)\Big] e^{ C Q(t)},
\]
where 
$$V^\xi(t)~=~\sum |\sigma_\alpha\xi_\alpha|,\qquad 
Q^\xi(t)\doteq \sum_{(\alpha,\beta) \in \mathcal A} |\sigma_\alpha \sigma_\beta|(|\xi_\alpha| + |\xi_\beta|), \qquad Q(t)\doteq \sum_{(\alpha,\beta)\in \mathcal A}|\sigma_\alpha\sigma_\beta|.
$$
Here the first sum ranges over all approaching fronts inside $\Gamma_{kj}$ at time $t$, while the 
other two sums range over all couples of approaching fronts, located inside $\Gamma_{kj}$.
In addition, $Q_{same}^\xi(t)$ denotes the amount of approaching waves with shifts, but restricted to
fronts of the same family:
\[
Q^\xi_{same}(t)\doteq \sum_{(\alpha,\beta) \in \A_{same}(t)} |\sigma_\alpha \sigma_\beta|(|\xi_\alpha| + |\xi_\beta|).\]
By (\ref{Linorm}), the speeds of any two fronts of the same family, located at $x_\alpha(t), x_\beta(t)$, satisfy
\bel{speed}\bigl|\dot x_\alpha(t)-\dot x_\beta(t)\bigr|~\leq~M_1 \sqrt{\ve_1}.\eeq
In view of (\ref{speed}), we thus define  $\A_{same}(t)$ as the set of all couples of wave fronts of the same
family (one of which is a shock), both contained in the trapezoid $\Gamma_{kj}$, and such that
\bel{disab}
\bigl| x_\alpha(t)-x_\beta(t)\bigr|~\leq~M_1 \sqrt{\ve_1} (t_{k+1}-t).
\eeq
Indeed, if their distance is larger than the right hand side of (\ref{disab}), the two fronts cannot interact during the 
time interval $[t, t_{k+1}]$.
\v
{\bf 2.} We now use the interaction estimates in Lemma~\ref{l:21}, recalling that by (\ref{Linorm}) the strength
of every wave-front is bounded by $|\sigma|=\O(1)\cdot \sqrt{\ve_1}$.
We claim that,
for a suitable choice of the constant $K$ and a sufficiently large $C$ (to be determined later), the change in the functional $ \sqrt \ve_1 ( V^\xi(t) + K Q^\xi(t)) e^{C  Q(t)}$
across any front interaction is bounded by
\[
\Delta \Big[ \sqrt \ve_1 \bigl( V^\xi(t) + K Q^\xi(t)\bigr) e^{CQ(t)}\Big]~\le ~-  e^{CQ(t-)}\frac{C}2 \sqrt{ \ve_1} |\sigma'\sigma''| V^\xi(t-).
\]
 Indeed, repeating the arguments  in \cite{BC95} one obtains:
\[
\bega{rl}
\Delta V^\xi &=~ \O(1)\cdot |\sigma'\sigma''| (|\xi'|+|\xi''|),\\[2mm]
\Delta Q^\xi&=\, -|\sigma'\sigma''|(|\xi'|+|\xi''|) + \O(1)\cdot|\sigma'\sigma''|(|\xi'|+|\xi''|) V(t-) \\[1mm]
&\qquad \quad + \O(1)\cdot |\sigma'\sigma''|(|\sigma'|+|\sigma''|) V^\xi(t-), \\[2mm]
 \Delta e^{CQ} &=\, -C|\sigma'\sigma''| e^{CQ}.
\enda
\]
By choosing $\delta_0$ sufficiently small, so that the total amount $V$ of waves inside the trapezoid $\Gamma_{kj}$ remains small,
we have
\[
\Delta Q^\xi~=\, -\frac12|\sigma'\sigma''|(|\xi'|+|\xi''|)+ \O(1) |\sigma'\sigma''|(|\sigma'|+|\sigma''|) V^\xi(t-).
\]
Therefore we can choose $K$ large enough so that 
\[
\Delta (V^\xi + KQ^\xi) \le -\frac{K}4|\sigma'\sigma''|(|\xi'|+|\xi''|) + \O(1)\cdot K |\sigma'\sigma''|(|\sigma'|+|\sigma''|) V^\xi(t-).
\]
Using the bound $|\sigma'|+|\sigma''| = \O(1)\cdot \sqrt \ve_1$, choosing $C$ sufficiently large depending on $K$ and the system (\ref{1}), we obtain
\[
\Delta  \big[(V^\xi + KQ^\xi) e^{CQ}\big]\le - e^{CQ} \frac{C}2|\sigma'\sigma''| V^\xi(t).
\]

Next, we estimate the change in $Q^\xi_{same}$. If two fronts of the same family interact at time $t$, then
\[\bega{rl} 
\Delta Q^\xi_{same}(t) &= ~- |\sigma'\sigma''|(|\xi'| + |\xi''|) + \O(1) \cdot|\sigma'\sigma''|(|\xi'| + |\xi''|) V_{same}(t-) 
\\[2mm]
&\qquad\qquad + \O(1) \cdot |\sigma'\sigma''|(|\sigma'| + |\sigma''|) V^\xi(t-),\enda
\]
where $V_{same}(t-) =  \O(1)\cdot \sqrt{ \ve_1}$ 
denotes the total strength of waves of the same family approaching the interacting fronts.
Since also $|\sigma'|+|\sigma''| = \O(1)\cdot \sqrt{ \ve_1}$, by choosing $\ve_1>0$  sufficiently small  and $C$ sufficiently large we obtain
\[
\Delta \Ups_{kj} (t)~ \le ~- \frac12 |\sigma'\sigma''|(|\xi'| + |\xi''|). 
\] 
Similarly,  if an interaction between fronts of different families occurs at time $t$, one has
\[
\Delta \Ups_{kj}(t)~ \le~ 0.
\] 
\v
{\bf 3.} Thanks to the previous estimates,
the total amount of interactions with shifts, among fronts of the same family contained inside  
the domain $\Gamma_{jk}$, can thus be bounded as
\bel{Qsa}
\Q^\xi_{same}(\Gamma_{kj})~ \le ~2\Ups_{kj}(t_k) ~=~ \O(1)\cdot \sqrt{\ve_1} \cdot V^\xi(t_k).
\eeq
By the analysis in Section~\ref{sec:3},  the total amount of wave strengths with shifts at time $t_k$ 
satisfies the the inductive estimates (\ref{shb2}).   In particular, for any $t_k\in [\tau,1]$ we have
\bel{Qs2}\sum_j V^\xi(t_k)~=~\O(1) \cdot (2L)^N |\sigma_0\xi_0|,\eeq
where $N = \O(1)\cdot |\ln\tau|$.
 
Summing over all trapezoids $\Gamma_{kj}$, from (\ref{Qsa}) and (\ref{Qs2})  we obtain the bound (\ref{HatQ}),
where the constant $C_\tau$ depends only on  $\tau$ and on the system (\ref{1}).
\endproof

\section{Well posedness for Temple class systems with $\L^\infty$ data}
\label{sec:4}
\begin{figure}[ht]
  \centerline{\hbox{\includegraphics[width=14cm]{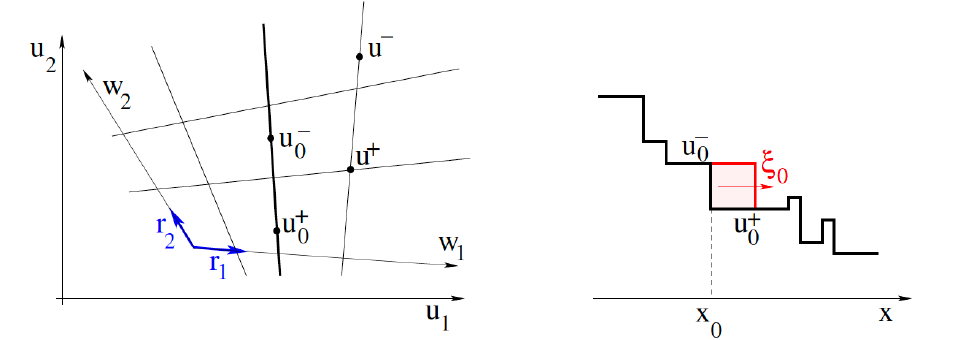}}}
  \caption{\small Left: right eigenvectors $r_1,r_2$ and Riemann invariants $w_1,w_2$ for a Temple class system.
  Right: a perturbation of the initial data, obtained by shifting the position of the jump initially located at $x_0$.}\label{f:ag83}
\end{figure}

\setcounter{equation}{0}
Let $(w_1, w_2)$ be a set of Riemann coordinates for the system (\ref{1}). 
Calling $u=u(w_1, w_2)$, this implies that the vectors 
\bel{r12}r_1 = {\partial u\over\partial w_1}\,,\qquad\qquad r_2 = {\partial u\over\partial w_2}\eeq
are eigenvectors of the Jacobian matrix $A(u)= Df(u)$.
We assume that both characteristic fields are genuinely nonlinear, and that the system
is of Temple class, namely all
integral curves of the vector fields $r_1$ and $r_2$ are straight lines
(see Fig.~\ref{f:ag83}, left).
In this setting, it is known that the system (\ref{1}) generates a globally Lipschitz semigroup, on any
positively invariant domain of the form
\bel{Dab}\bega{l}
\D~\doteq~\Big\{ u\in  \L^1(\R\,;\R^2)\,;~~~u(x)=u(w_1(x), w_2(x)\bigr)\\[2mm]
\qquad\qquad\qquad\qquad~~\hbox{with}~~w_1(x)\in [a_1, b_1], 
~w_2(x)\in [a_2, b_2]~~\forall x\in\R\Big\}.\enda\eeq
Notice that this is a
domain of bounded integrable functions, which may well have infinite total variation.

\begin{theorem}\label{t:31} {\bf (Semigroup generated by a Temple class system).} 
In the above setting, there exists a semigroup of entropy weak solutions to the system (\ref{1}),
which is globally Lipschitz continuous w.r.t.~the initial data. Namely, there exists a constant $L^*$ such that
\bel{glolip} \bigl\| S_t \bar u- S_t\bar v\bigr\|_{\L^1}~\leq~L^*\, \|\bar u-\bar v\|_{\L^1}\,,\eeq
for every $t\geq 0$ and every initial data $\bar u,\bar v\in \D$.
\end{theorem}

{\bf Proof.} 
This result was first proved in \cite{BG2}. We give here a simpler argument, that leads to better estimates on the 
Lipschitz constant $L^*$.
The property of Temple class
implies that 
\begi
\item[(i)] Rarefaction curves and shock curves coincide, since they are all straight lines.
\item[(ii)] The strength of wave fronts, measured in terms of the change of the corresponding Riemann
invariant, does not change across interactions.
More precisely,  in the interaction of fronts of different families shown in Fig.~\ref{f:ag51} left, 
one has
\bel{Td}\sigma_1\,=\, \sigma',\qquad\sigma_2 \,=\,\sigma''.\eeq
In the interaction of fronts of the same family (say, the second one) shown in Fig.~\ref{f:ag51} right, 
one has
\bel{Ts}\sigma_1\,=\, 0,\qquad\sigma_2 \,=\,\sigma'+\sigma''.\eeq
\endi
Following the approach introduced in \cite{BC95}, 
to estimate the Lipschitz constant of the semigroup it suffices to prove that the same Lipschitz continuous
dependence holds for all front tracking approximations. 

We thus fix an integer $\nu\geq 1$ and consider a piecewise constant front tracking solution 
constructed as in \cite{BC95}. 

Consider one single front, say located at $x_0$ at time $t=0$.   
Call $\sigma_0$ the size of this front, and consider a perturbation where this front is initially shifted at rate $\xi_0$.
In turn, for $t>0$, this will produce shifts in all
wave fronts in its domain of dependency.  If we can find a constant $L$ such that
the sizes $\sigma_\alpha(t)$ and the shifts $\xi_\alpha(t)$  of all fronts at any time $t>0$
satisfy
\bel{Lip8}\sum_\alpha \bigl|\sigma_\alpha(t)\bigr| \, \bigl|\xi_\alpha(t)\bigr|~\leq~L\cdot |\sigma_0|\, |\xi_0|\,,\eeq
this will establish the Lipschitz continuity of the semigroup w.r.t.~the initial data, in the $\L^1$ distance.

 To fix ideas, let a time interval $[0,T]$ be given 
and consider a small shift $\xi_0$ of a front of the second family, initially
located at $x_0$, with size $\sigma_0$.  
The key estimate (\ref{Lip8}) will be established by analyzing four possible configurations.
In each case, we shall  compute exactly the shifts $\xi$ of all the fronts at time $t=T$, depending on their position with respect to the front shifted at time $t=0$. The expression will depend on the left and right states $u^-,u^+$ across the front at time $T$, the left and right states $u^-_0,u^+_0$ across the front which is shifted at the initial time $t=0$, and the initial shift $\xi_0$.

\begin{figure}[ht]
  \centerline{\hbox{\includegraphics[width=16cm]{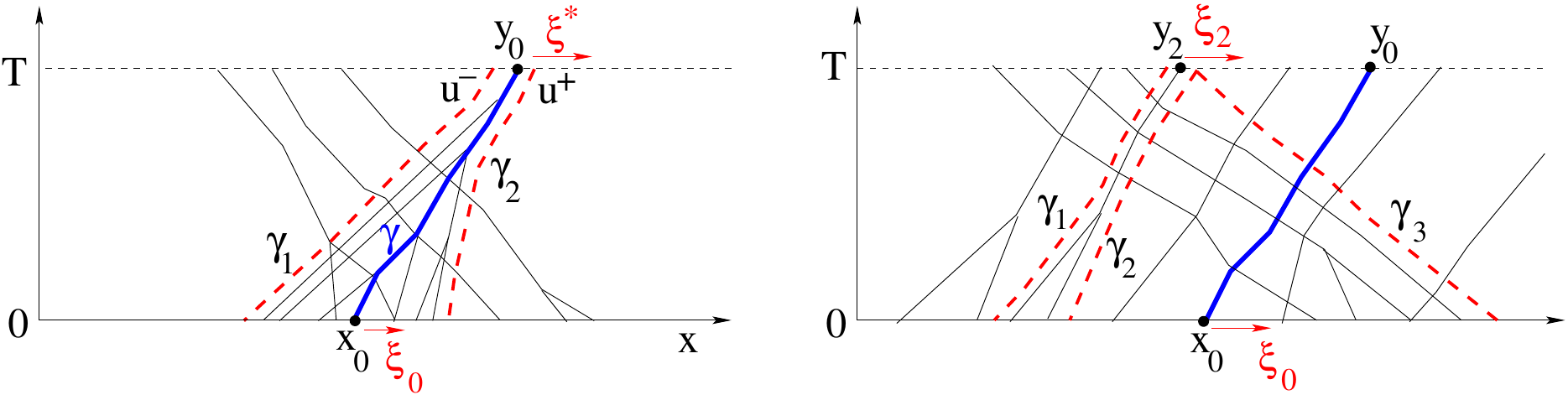}}}
  \caption{\small Left: computing the shift $\xi^*$ at time $T$ of the
  2-front that initiates at $x_0$.   This is achieved by using the divergence theorem on the region bounded by the curves $\gamma_1, \gamma_2$.
  Right:  computing the shift $\xi_2$ of the 2-front located at $y_2$
  at time $T$.
   }\label{f:ag82}
\end{figure}

As in (\ref{r12}), we denote by $r_1(u),r_2(u)$ a basis of right eigenvectors of $Df(u)$ and by $\ell_1(u),\ell_2(u)$ the dual basis, so that
\bel{dbasis} \ell_i(u)\cdot r_j(u)~=~\left\{ \bega{rl} 1\quad &\hbox{if} ~~i=j,\cr
0\quad &\hbox{if} ~~i\not=j.
\enda\right. \eeq

{\bf CASE 1.} We start by calculating the shift at the terminal time $T$ of the 2-front $\gamma$ starting at $x_0$.  As shown in Fig.~\ref{f:ag82}, left, let $y_0$  be the terminal position of this front. Consider a curve $\gamma_1$ slightly on the left of the minimal backward characteristic from $y_0$ and a curve $\gamma_2$ slightly on the right of the maximal backward characteristic from $y_0$.

We denote by $\Delta f\bigr|_{\gamma_1}$ the difference in the total flux across $\gamma_1$ 
between the solution $u^\xi$ with the shifted initial data  and the original solution $u$. 
Similarly, we denote by  $\Delta f\bigr|_{\gamma_2}$ the difference in the total flux across $\gamma_2$, 
between the two solutions.   Notice that, since the perturbation does not affect points to the right 
of the curve $\gamma$, we trivially have $\Delta f\bigr|_{\gamma_2}=0$.

Since both solutions $u^\xi$, $u$ satisfy the conservation equations (\ref{1}), applying the divergence theorem  on the region 
$$\Omega ~\doteq~ \bigl\{(t,x) \in [0,T]\times \R \,;~ \gamma_1 (t)< x<\gamma_2(t) \bigr\},$$ 
one obtains
\bel{e:cons1}
\Delta f\bigr|_{\gamma_1} + (u^-_0 - u^+_0) \xi_0 ~=~  (u^--u^+) \xi^*.
\eeq
Here $\xi^*$ is the shift of the front $\gamma$ at the terminal time $T$, and $u^-, u^+$ are the 
terminal left and right states, respectively.

We now observe that all fronts crossing $\gamma_1$ belong to the first family. Hence  the second Riemann invariant $w_2$ is constantly equal to $w_2 (u^-)$ along $\gamma_1$. 
The difference in flux is thus parallel to the right eigenvector:
\bel{parl1}
\Delta f\bigr|_{\gamma_1}~ \parallel  ~ r_1(u^-).
\eeq
Multiplying  both sides of \eqref{e:cons1} by the left eigenvector $\ell_2(u^-)$, in view of 
(\ref{dbasis})  one obtains
\[
\ell_2(u^-) \cdot (u^-_0 - u^+_0) \xi_0 ~=~ \ell_2(u^-) \cdot (u^--u^+) \xi^*,
\]
hence
\bel{shift1} \xi^*~=~{\ell_2(u^-) \cdot (u^-_0 - u^+_0)\over \ell_2(u^-) \cdot (u^--u^+)}\,\xi_0\,.
\eeq
Toward an approximation estimate, we write
$$\sigma^* ~=~w_2(u^+)-w_2(u^-)~=~\int_{w_2(u^-)}^{w_2(u^+)}~1\, ds~=~
\int_{w_2(u^-)}^{w_2(u^+)}~\ell_2\bigl( u(w_1,s)\bigr)\cdot r_2\bigl( u(w_1,s)\bigr)\, ds,$$
where $w_1\doteq w_1(u^+)=w_1(u^-)$.
On the other hand,
$$u^+-u^-~=~\int_{w_2(u^-)}^{w_2(u^+)}r_2\bigl( u(w_1,s)\bigr)\, ds,$$
hence
\bel{s10}\bega{l} \bigl| \sigma^* - \ell(u^-)\cdot (u^+-u^-)\bigr| ~\leq~|u^+-u^-|\cdot \max_{s} \Big| \ell_2(u(w_1,s)) - \ell_2(u^-)\Big|\\[3mm]
\qquad=~\O(1)\cdot \|u\|_{\L^\infty} \cdot |u^+-u^-|~=~\O(1)\cdot \|u\|_{\L^\infty}\cdot \sigma^*.\enda\eeq
A similar argument yields
\bel{s11}\bigl| \sigma_0 - \ell(u^-)\cdot (u_0^+-u_0^-)\bigr| ~=~\O(1)\cdot \|u\|_{\L^\infty} \cdot \sigma_0\,.\eeq
Combining the two previous bounds we conclude
\bel{sb1}
|\sigma^* \xi^* - \sigma_0 \xi_0| ~= ~ \O(1)\cdot \|u\|_{\L^\infty} \cdot \bigl(| \sigma^*\xi^*| +|\sigma_0 \xi_0|\bigr)~=~  \O(1)\cdot \|u\|_{\L^\infty}  \cdot |\sigma_0 \xi_0|\,.
\eeq

{\bf CASE 2.}  
 As shown in Fig.~\ref{f:ag82}, right, we now compute the shift $\xi_2$ of a 2-front, which at the terminal time $T$ is located at  the point $y_2 < y_0$.

Let $\gamma_1$ be a path slightly to the left of the minimal backward 2-front ending at $y_2$ and 
let $\gamma_2$ be a path slightly to the right of the maximal backward 2-front ending at $y_2$. 
Moreover, let $\gamma_3$ be a path slightly to the right of the first maximal backward 1-front which at time $T$ ends  to the left of $y_2$.
We only need to consider the case where $\gamma_3(0)>  x_0$, otherwise the shift in $x_0$
does not affect the solution in a neighborhood of $y_2$.

As before, we denote by $\Delta f\bigr|_{\gamma_i}$, $i=1,2,3$, the difference in the total flux across $\gamma_i$ 
between the solution $u^\xi$ with the shifted initial data  and the original solution $u$. 

Since both solutions satisfy the conservation equations (\ref{1}), applying the divergence theorem  on the two regions
\[
\Omega_1 = \bigl\{(t,x) \in [0,T]\times \R \,;~~ \gamma_1(t)<x<\gamma_2(t) \bigr\}, \qquad 
\Omega_2 = \bigl\{(t,x) \in [0,T]\times \R\,;~~ \gamma_2(t) <x < \gamma_3(t)\bigr\},
\]
we obtain
\bel{s45}\Delta f\bigr|_{\gamma_1}-\Delta f\bigr|_{\gamma_2}~=~ (u^--u^+) \xi_2\,,
\qquad\quad\Delta f\bigr|_{\gamma_3}\,=\,  (u^-_0 -u^+_0)\xi_0 + \Delta f\bigr|_{\gamma_2}\,.\eeq
Here $u^-,u^+$ denote the left and right states across the 2-front at $y_2$, while $\xi_2$ is its shift.
We observe that the following parallelism relations hold:
\bel{parl2}
\Delta f\bigr|_{\gamma_1} \parallel r_1 (u^-), \qquad \quad\Delta f\bigr|_{\gamma_2}\parallel r_1(u^+), \qquad \quad\Delta f\bigr|_{\gamma_3}\parallel r_2(u^+).
\eeq
Indeed, the curves $\gamma_1$ and $\gamma_2$ are crossed only by 1-fronts, while $\gamma_3$
is only crossed by 2-fronts.
Taking the product with $\ell_1(u^+)$, from the second equation in (\ref{s45}) it follows 
\bel{s46}
-\ell_1(u^+)\cdot \Delta f\bigr|_{\gamma_2}~=~\ell_1(u^+)\cdot (u^-_0 -u^+_0)\xi_0 \,.\eeq
In view of the second relation in (\ref{parl2}), this implies
\bel{Df2}\Delta f\bigr|_{\gamma_2}~=~-\Big(\ell_1(u^+)\cdot (u^-_0 -u^+_0)\Big)\xi_0\, r_1(u^+).\eeq
In turn, multiplying the first equation in (\ref{s45}) by $\ell_2(u^-)$ we obtain
\bel{s47}
\ell_2(u^-) \cdot (u^--u^+) \xi_2~=~-\ell_2(u^-) \cdot \Delta f\bigr|_{\gamma_2}~=~
\Big( \ell_2(u^-) \cdot r_1(u^+)\Big)\Big(\ell_1(u^+)\cdot (u^-_0 -u^+_0)\Big)\xi_0\,.\eeq
Therefore
\bel{shift2}
\xi_2~=~{
\Big( \ell_2(u^-) \cdot r_1(u^+)\Big)\Big(\ell_1(u^+)\cdot (u^-_0 -u^+_0)\Big)
\over \ell_2(u^-) \cdot (u^--u^+)}\, \xi_0\,.\eeq

Since 
$$\bigl|\ell_2(u^-) \cdot r_1(u^+)\bigr| ~=~ \bigl|(\ell_2(u^-) - \ell_2(u^+))\cdot r_1(u^+)\bigr| ~=~ \O(1)\cdot \sigma,$$ 
from (\ref{shift2}) we obtain the estimate
\bel{sb2}
|\sigma \xi_2| ~=~ \O(1)\cdot |\sigma| |\sigma_0 \xi_0|.
\eeq

\begin{figure}[ht]
  \centerline{\hbox{\includegraphics[width=16cm]{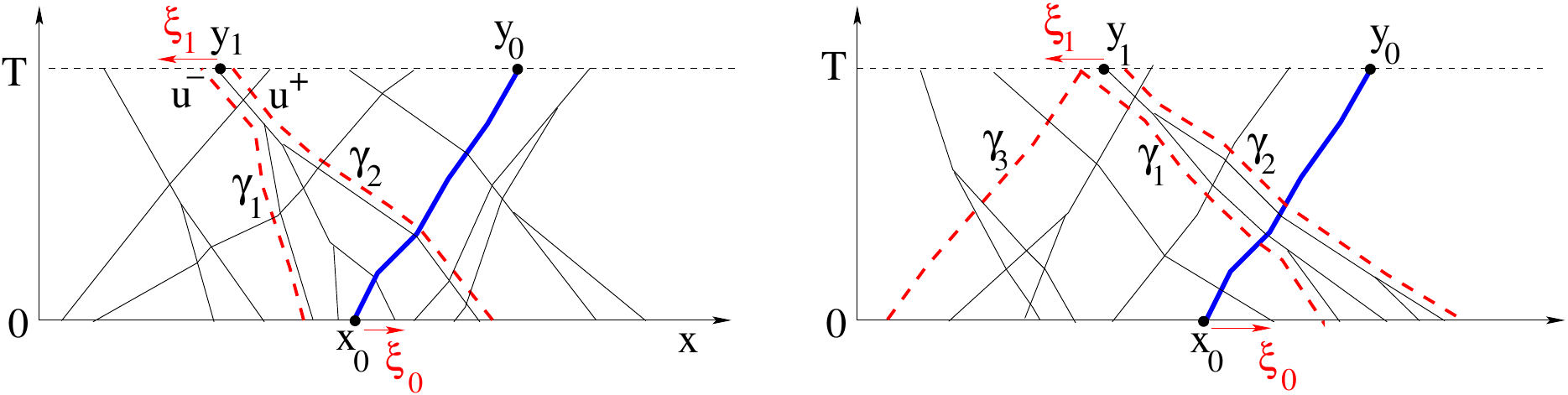}}}
  \caption{\small Computing the shift $\xi_1$ of the  1-front located at $y_1$ at time $T$.
  Left: the construction used in Case 3a. Right: the construction used in Case 3b.}\label{f:ag84}
\end{figure}

{\bf CASE 3.} We now consider the shift $\xi_1$ of a 1-front,   which at the terminal time $T$ is located at a point 
 $y_1 < y_0$.  As shown in Figure~\ref{f:ag84}, let $\gamma_1$ be a path slightly to the left of the minimal backward 1-front ending at $y_1$ and 
let $\gamma_2$ be a path slightly to the right of the maximal backward 1-front ending at $y_1$. 
We only study the case where $\gamma_1(0)<  x_0$, because otherwise the shift in $x_0$
does not affect the solution in a neighborhood of $y_1$.
Two sub-cases need to be considered.

{\bf Case 3a.} Assume $\gamma_1(0)< x_0<\gamma_2(0)$ (see Fig.~\ref{f:ag84}, left). 
Notice that this situation can occur for at most one terminal point $y_1$. 
In this case the initial shift in the 2-front at $x_0$
does not affect the flux through $\gamma_1$, hence $\Delta f\Big|_{\gamma_1}=0$.
Since both solutions $u^\xi$, $u$ satisfy the conservation equations (\ref{1}), applying the divergence theorem  on the region 
 \[
 \Omega \,= \, \bigl\{(t,x) \in [0,T]\times \R \,;~ \gamma_1 (t)< x<\gamma_2(t) \bigr\}
 \]
 one obtains
 \bel{s48}
 (u^-_0 -u^+_0)\xi_0 ~=~ \Delta f\Big|_{\gamma_2} + (u^--u^+) \xi_1.
 \eeq
 Here $u^-,u^+$ denote the left and right states across the 1-front at $y_1$, while $\xi_1$ is its shift.
 Since $\Delta f\Big|_{\gamma_2}$ is parallel to $ r_2(u^+)$, multiplying both sides of (\ref{s48}) by 
 $\ell_1(u^+)$ we obtain
 \bel{sh4}
 \ell_1(u^+)\cdot (u^-_0 -u^+_0)\xi_0 ~=~ \ell_1(u^+) (u^--u^+) \xi_1\,.
 \eeq
Therefore
 \bel{shift3a}
 \xi_1 ~=~ \frac{\ell_1(u^+) \cdot (u^-_0-u^+_0)}{\ell_1(u^+)\cdot  (u^--u^+)}\xi_0 \,.
 \eeq
 To estimate the size of this term, we observe that $u_0^--u_0^+$ is parallel to $r_2(u_0^-)$.
 Hence
 $$\bega{l}
 \ell_1(u^+) \cdot  (u^-_0-u^+_0) ~=~ \bigl( \ell_1(u^+) - \ell_1(u_0^-)\bigr)  \cdot  (u^-_0-u^+_0)\\[3mm]
 \quad =~\O(1)\cdot (u^+- u_0^-)  (u^-_0-u^+_0)~=~\O(1)\cdot \|u\|_{\L^\infty}\cdot \sigma_0\,.
 \enda$$
 Calling $\sigma=\O(1)\cdot |u^+-u^-|$ the strength of the jump at $y_1$, from (\ref{sh4}) we deduce
\bel{sb3a}
|\sigma \xi_1| ~=~ \O(1)\cdot \|u\|_{\L^\infty}\cdot |\sigma_0 \xi_0|.
\eeq
\v
{\bf Case 3b.} We now assume $\gamma_1(0)> x_0$ (see Fig.~\ref{f:ag84}, right).
Define the paths  $\gamma_1$ and $\gamma_2$  as in Case 3a, and let $\gamma_3$ be a path
slightly to the left of the minimal backward 2-characteristic which at time $T$ reaches $y_1$. 
Since both solutions $u^\xi$, $u$ satisfy the conservation equations (\ref{1}), applying the divergence theorem  on the two regions
\[\bega{rl}
\Omega_1&=~ \bigl\{ (t,x) \in [0,T]\times \R \,;~ \gamma_3(t)<x<\gamma_1(t)\bigr\}, \\[2mm]
\Omega_2&=~ \bigl\{ (t,x) \in [0,T]\times \R \,;~ \gamma_1(t)<x<\gamma_2(t)\bigr\},\enda
\]
we obtain
\bel{shi7}
 \Delta f\Big|_{\gamma_1}~ = ~ \Delta f\Big|_{\gamma_3} + (u^-_0-u^+_0)\xi_0\,,
 \qquad\quad 
 \Delta f\Big|_{\gamma_1}~ =  \Delta f\Big|_{\gamma_2} +  (u^--u^+) \xi_1 \,.
\eeq
As before, $u^-, u^+$ denote the left and right states across the 1-front  at $y_1$, while $\xi_1$
is its shift.
Using the relations
\[
 \Delta f\Big|_{\gamma_1} \parallel r_2(u^-), \qquad \Delta f\Big|_{\gamma_2} \parallel r_2(u^+), \qquad  \Delta f\Big|_{\gamma_3}\parallel r_1(u^+),
\]
multiplying the first equation in (\ref{shi7}) by $\ell_2(u^+)$ we find
\[\ell_2(u^+)\cdot  \Delta f\Big|_{\gamma_1}~ = ~ \ell_2(u^+)\cdot (u^-_0-u^+_0)\xi_0\,.\]
Hence
\bel{so1}
\Delta f_{\gamma_1}~ =~ \frac{\ell_2(u^+)\cdot (u_0^- - u_0^+) \xi_0}{\ell_2(u^+)\cdot r_2(u^-)} r_2(u^-).\eeq
Multiplying the second equation in (\ref{shi7}) by $\ell_1(u^+)$ and then using (\ref{so1})
we obtain
\[ \ell_1(u^+)\cdot (u^--u^+) \xi_1~ =~ \ell_1(u^+)\cdot \Delta f\Big|_{\gamma_1}~=~
\Big(\ell_1(u^+) \cdot r_2(u^-)\Big)\frac{\ell_2(u^+) \cdot (u_0^- - u_0^+)}{\ell_2(u^+) \cdot r_2(u^-)}\xi_0\,,
\]
\bel{shift3b}
\xi_1~=~{\ell_1(u^+) \cdot r_2(u^-)\over  \ell_1(u^+)\cdot (u^--u^+) }\cdot \frac{\ell_2(u^+) \cdot (u_0^- - u_0^+)}{\ell_2(u^+) \cdot r_2(u^-)}\,\xi_0\,,
\eeq
 Calling $\sigma=\O(1)\cdot |u^+-u^-|$ the strength of the jump at $y_1$, 
and observing that 
$$ \ell_1(u^+) \cdot r_2(u^-)~=~\bigl(  \ell_1(u^+) -\ell_1(u^-)\bigr)\cdot r_2(u^-) ~=~ \O(1)\cdot \sigma,$$ 
by (\ref{shift3b}) we conclude
\bel{sb3b}
|\sigma \xi_1| = \O (1)\cdot |\sigma|\, |\sigma_0 \xi_0|.
\eeq
\v
We are now ready to prove the key estimate (\ref{Lip8}).  
Consider a front tracking solution $u$, and a perturbed solution $u^\xi$ obtained by 
shifting a front initially located at $x_0$.  Let 
$\sigma_0$ and $\xi_0$ be respectively the size and shift of this front.  
Since we are assuming that all wave speeds are $\leq 1$ in absolute value, at a later
time $T>0$ the only fronts that can be shifted in the perturbed solution will be the ones located at
points $x_\alpha\in [x_0-T, x_0+T]$.   
In view of (\ref{sb1}), (\ref{sb2}), (\ref{sb3a}), (\ref{sb3b}), we obtain
\bel{sh5}\sum_{x_\alpha\in  [x_0-T,\, x_0+T]} |\sigma_\alpha\xi_\alpha|~=~
\O(1)\cdot\bigl( 1 + \|u\|_{\L^\infty}\bigr) |\sigma_0\xi_0| + 
\O(1)\cdot
\sum_{x_\alpha\in  [x_0-T,\, x_0+T]} |\sigma_\alpha|\cdot |\sigma_0 \xi_0|.\eeq
Because of genuine nonlinearity, the total strength of waves decays in time.
In particular, for some constant $C_0$ there holds
\bel{e:tv-gnl}\sum_{x_\alpha\in  [x_0-T,\, x_0+T]}  |\sigma_\alpha|~\leq~C_0\,,\eeq
for every $x_0$ and every $T>0$.
Inserting this bound in (\ref{sh5}) we conclude
$$\sum_{x_\alpha\in  [x_0-T,\, x_0+T]} |\sigma_\alpha\xi_\alpha|~=~\O(1)\cdot |\sigma_0\xi_0|.$$
This yields (\ref{Lip8}), completing the proof.
\endproof

\section{Comparison with a Temple class system} 
\label{sec:5}
\setcounter{equation}{0}
Starting with the general $2\times 2$ system (\ref{1}), consider the 
diagonalized system (\ref{diag}), written in Riemann invariants.
We wish to approximate this system by a Temple class system, say 
\bel{TCS}
u_t + \Tilde f(u)_x~=~0,\eeq
in a neighborhood of the origin.  The following lemma shows that, by a suitable choice of the 
Temple system,  the first order derivatives of the eigenvalues of the two systems coincide at the origin.

\begin{lemma}\label{l:41} Consider a strictly hyperbolic $2\times 2$ system (\ref{1}), written in Riemann coordinates as in (\ref{diag}).
Then there exists a $2\times 2$ Temple class system, whose integral curves of eigenvectors are straight lines,  which in Riemann
invariants takes the form
\bel{Tedi}\left\{ \bega{rl}w_{1,t}+\mu_1(w_1, w_2) w_{1,x}&=~0,\\[1mm]
w_{2,t}+\mu_2(w_1,w_2) w_{1,x}&=~0,\enda\right.
\eeq
with 
\bel{tanTe}\mu_i(0,0)\,=\,\lambda_i(0,0),\quad
\qquad c_{ij}~\doteq~{\partial \over\partial w_j}\mu_i(0,0) ~=~ {\partial \over\partial w_j}\lambda_i(0,0),
\qquad i,j=1,2.\eeq
\end{lemma}

\begin{figure}[ht]
  \centerline{\hbox{\includegraphics[width=12cm]{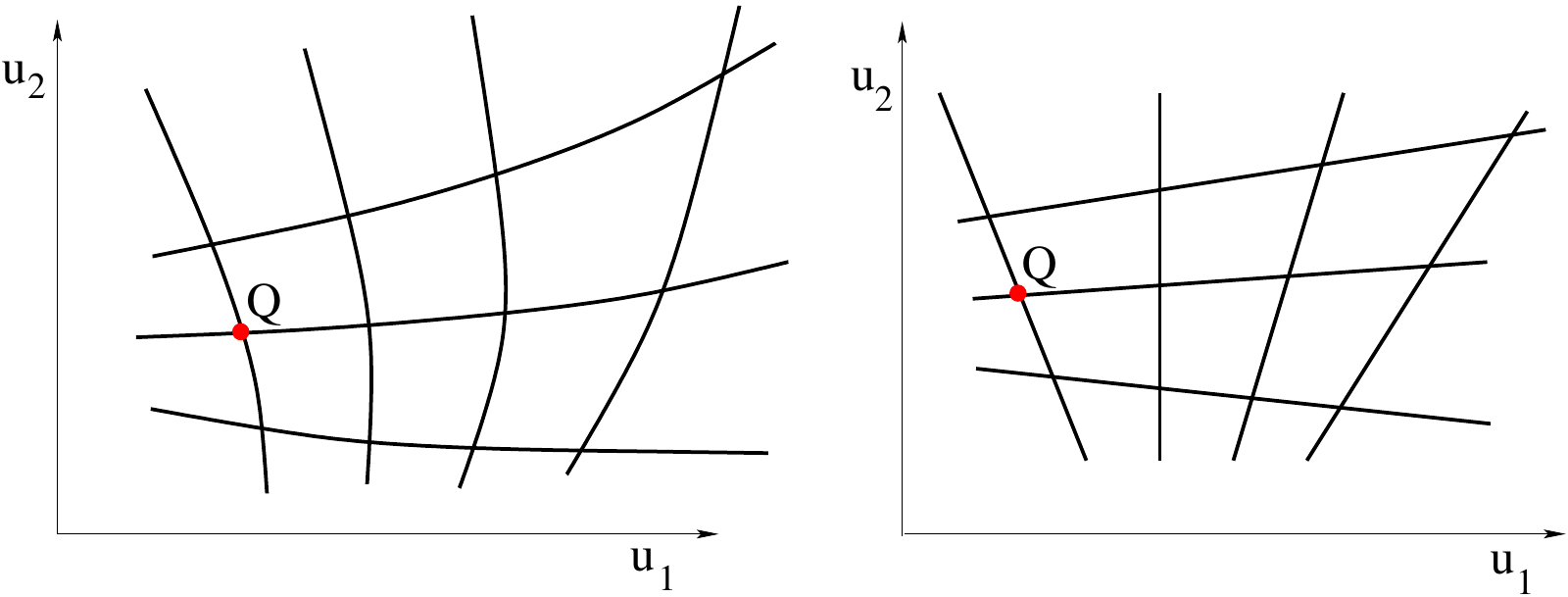}}}
  \caption{\small Left:  Riemann coordinates for the original system (\ref{1}).
  Right: a Temple class system which is ``tangent" to the original one at 
  the point $U_0$.}\label{f:ag57}
\end{figure}
{\bf Proof.} 
Consider a Temple class system, with  Riemann coordinates $(w_1, w_2)$ and eigenvectors
$r_1,r_2$ as in (\ref{r12}).   In this proof, it will be convenient to work with 
normalized eigenvectors $\bfr_1,\bfr_2$ having unit length.
Since integral curves of eigenvectors are straight lines, these normalized eigenvectors satisfy
$$\bfr_1\,=\,\bfr_1(w_2),\qquad \bfr_2\,=\,\bfr_2(w_1)$$
In other words, the first eigenvector does not depend on $w_1$ and the second eigenvector does not depend on $w_2$. By construction we have
$${\partial \Tilde f\over\partial w_1}(w_1, w_2) ~=~\alpha(w_1, w_2) \bfr_1(w_2),
\qquad\qquad {\partial \Tilde f\over\partial w_2}(w_1, w_2) ~=~\beta(w_1, w_2) \bfr_2(w_1)$$
for some scalar functions $\alpha,\beta$.
This implies
\bel{2der}
{\partial^2 \Tilde f\over\partial w_1\partial w_2}~=~\alpha_{w_2} \bfr_1 + \alpha \bfr_{1, w_2}
~=~\beta_{w_1} \bfr_2 + \beta \bfr_{2, w_1}\,.\eeq
Notice that (\ref{2der}) is a vector equation. It can be reduced to two scalar equations
letting $\{l_1,l_2\}$ be the dual basis of $\{\bfr_1, \bfr_2\}$.    Taking the inner products of 
(\ref{2der}) with $l_1$ and $l_2$ respectively, one obtains
\bel{absy}
\left\{
\bega{rl} \alpha_{w_2} &=~\bigl( l_1\cdot \bfr_{2, w_1}\bigr) \beta - 
\bigl( l_1\cdot \bfr_{1, w_2}\bigr)\alpha ,\\[2mm]
\beta_{w_1}&=~\bigl( l_2\cdot \bfr_{1, w_2}\bigr)\alpha  -  \bigl( l_2\cdot \bfr_{2, w_1}\bigr)\beta.
\enda\right.
\eeq
The system (\ref{absy}) can be solved in a neighborhood of the origin, with 
any smooth boundary data of the form
\bel{bda}
\alpha(w_1, 0)~=~\bar\alpha(w_1),\qquad\qquad \beta(0, w_2)~=~\bar\beta(w_2).\eeq
This already implies that the coefficients $\alpha_{w_1}$ and $\beta_{w_2}$
at the origin can be chosen arbitrarily.
In turn, since the values $l_1\cdot \bfr_{2,w_1}$ and $l_2\cdot \bfr_{1,w_2}$ can 
be given arbitrarily   (by suitably choosing the fields of unit vectors $\bfr_2(w_1)$
and $\bfr_1(w_2)$), this implies that also the 
derivatives
$\alpha_{w_2}$ and $\beta_{w_1}$ at the origin can be assigned arbitrarily.
\endproof

Next, consider the interaction of two fronts, with strengths $\sigma', \sigma''$ and 
shifts $\xi', \xi''$, as in Fig.~\ref{f:ag51}. We shall estimate the difference between:
\begi
\item[(i)] the strengths  $\sigma_1, \sigma_2$ and 
shifts $\xi_1, \xi_2$ of the outgoing fronts, for the original system, where the left state
is $u^-=u(w_1^-, w_2^-)$ as in (\ref{lrstate}).
\item[(ii)] the strengths  $\Tilde \sigma_1, \Tilde \sigma_2$ and 
shifts $\Tilde \xi_1, \Tilde \xi_2$ of the outgoing fronts, for the corresponding solution of the Temple system, where the left state
is $\Tilde u^-=u(\Tilde w_1^-, \Tilde w_2^-)$.
\endi

\begin{lemma}\label{l:52}
Consider the interaction of two fronts, as in Fig.~\ref{f:ag51}.
Assume that the left states satisfy $|\Tilde u^--u^-|\leq C_1 \ve$ and that the incoming fronts have strength
$|\sigma'|, |\sigma''|\leq C_1 \ve$.\\
Then the difference in strengths and shifts of the outgoing fronts, 
between the original system (\ref{1}) and the Temple system (\ref{TCS}) satisfying (\ref{tanTe})
satisfy the following bounds.

CASE 1: If the two incoming fronts belong to different families, then
\bel{bT1} |\Tilde \sigma_1 - \sigma_1|+ |\Tilde \sigma_2-\sigma_2|~=~\O(1)\cdot \ve^2\,|\sigma'\sigma''|\,,
\eeq
\bel{bT2}|\Tilde \xi_1-\xi_1|~=~\O(1) \cdot \ve |\sigma''|\, \bigl(|\xi'|+|\xi''|\bigr), 
\qquad \quad |\Tilde \xi_2-\xi_2|~=~\O(1) \cdot \ve |\sigma'|\, \bigl(|\xi'|+|\xi''|\bigr).
\eeq

CASE 2: If the two incoming fronts belong to the same  family (say, the first one) then
\bel{bT3} |\Tilde \sigma_1 - \sigma_1|+ |\Tilde \sigma_2-\sigma_2|~=~\O(1)\cdot \ve\,|\sigma'\sigma''|\,,
\eeq
\bel{bT43}\Tilde \xi_1\Tilde \sigma_1  - \xi_1\sigma_1 ~=~\O(1)\cdot\,\bigl( |\xi'|+|\xi''|\bigr) |\sigma'\sigma''|,
\eeq
\bel{bT5}|\sigma_2 \xi_2|~=~\O(1) \cdot |\sigma'\sigma''| \bigl(|\xi'|+|\xi''|\bigr).\eeq
\end{lemma}

{\bf Proof.}
In Case 1, the bound (\ref{bT1}) is an immediate consequence of (\ref{ies1}),
because for a Temple system one has $\Tilde\sigma_1=\sigma'$ and $\Tilde\sigma_2=\sigma''$.

Concerning the shifts of the outgoing fronts, we recall that for the original system these shifts were 
computed at (\ref{shi1})-(\ref{shi2}). 
On the other hand, for the Temple system by (\ref{tanTe}) it follows
\bel{sf1} \Tilde \xi_1~=~\xi' + {c_{12} \sigma'' (\xi'-\xi'')\over \Tilde \lambda''-\Tilde\lambda'}
+\O(1) \cdot \ve |\sigma''|\, \bigl(|\xi'|+|\xi''|\bigr),\eeq
where $ \Tilde \lambda'$, $\Tilde\lambda''$ are the speeds of the incoming fronts.
By (\ref{tanTe}), the differences between these speeds are
\bel{TLL} |\Tilde \lambda'-\lambda'|\,=\,\O(1)\cdot \ve,
\qquad\qquad |\Tilde\lambda''-\lambda''|\,=\,\O(1)\cdot \ve .\eeq
Since by strict hyperbolicity the difference $\lambda''-\lambda'$ is bounded away from zero, 
this yields
\bel{sf2} \Tilde \xi_1~=~\xi' + {c_{12}\, \sigma'' (\xi'-\xi'')\over  \lambda''-\lambda'}
+\O(1) \cdot \ve |\sigma''|\, \bigl(|\xi'|+|\xi''|\bigr),\eeq
Combining (\ref{sf2}) with (\ref{shi1}) we conclude
$$|\Tilde \xi_1-\xi_1|~=~\O(1) \cdot \ve |\sigma''|\, \bigl(|\xi'|+|\xi''|\bigr).$$
An entire similar estimate holds for the difference $|\Tilde \xi_2-\xi_2|$. This yields
(\ref{bT2}). 
\v
Next, consider Case 2 where the interacting  fronts are of the same family, say, the first one.
The bound (\ref{bT2}) is an immediate consequence of (\ref{ies3}),
because for the Temple system one has $\Tilde\sigma_1=\sigma'+\sigma''$ and $\Tilde\sigma_2=0$.
In addition, from  (\ref{ies3})-(\ref{ies4}) it follows
\bel{bT41} 
|\sigma_2 \xi_2|~=~\O(1) \cdot |\sigma'\sigma''| \bigl(|\xi'|+|\xi''|\bigr).\eeq
The same argument used at (\ref{ospe})-(\ref{ies71}) can be applied to the Temple system as well, 
thus obtaining
\bel{bT42}\left|\Tilde  \xi_1-{\xi'\sigma'+\xi''\sigma''\over \sigma'+\sigma''} \right|~=~\O(1)\cdot \bigl( |\xi'|+|\xi''|\bigr) {|\sigma'\sigma''| \over |\sigma'| + |\sigma''|},\eeq
\bel{bT6} 
\Tilde \xi_1\Tilde \sigma_1  - \xi_1\sigma_1 ~=~\O(1)\cdot\bigl( |\xi'|+|\xi''|\bigr) |\sigma'\sigma''|\,.\eeq

\section{The Lipschitz constant of the semigroup on $\D^{(\tau)}$} 
\label{sec:6}

In this section we estimate the Lipschitz constant of the semigroup generated
by (\ref{1}), on  the domain $\D^{(\tau)}$ at (\ref{Dtau}).   
As explained in the Introduction, to achieve such an estimate one can
consider a front tracking approximation $u:[\tau, 1]\times \R\mapsto \R^2$, 
and shift one of the fronts at time $t=0$, say located at $x_0$, with strength $\sigma_0$ and  shift rate $\xi_0$.   By subsequent wave front  interactions, this will produce
shifts in many other fronts.   The goal is to estimate the quantity
\bel{mom1}\sum_\alpha \hbox{[strength] $\times$ [shift]}~\doteq~\sum_\alpha \bigl|\sigma_\alpha(t)
\xi_\alpha(t)\bigr|\,,\eeq 
where the sum ranges over all wave-fronts of $u(t,\cdot)$,  at  any time $t\in [\tau,1]$, see Fig.~\ref{f:ag68}.

For a Temple class system, the proof of Theorem~\ref{t:31} shows that there exists a constant $L^*$ such that 
the bound (\ref{Lip8}) holds for every $t\geq 0$ and every piecewise constant initial data
in the domain $\D$ at (\ref{Dab}).  

To prove a similar bound for the original $2\times 2$ system (\ref{1}), the key idea is that, as long as the
norm $\bigl\|u(t,\cdot)\bigr\|_{\L^\infty}$ remains small, the system (\ref{1}) can be 
closely  approximated by a Temple class system.  As a consequence, for any given $\tau>0$,
by choosing $\ve_1>0$ sufficiently small, on the domain $\D^{(\tau)}$ at (\ref{Dtau})  we will achieve a 
Lipschitz constant $L(\tau)\leq 2L^*$.

\begin{lemma}\label{l:61} Let (\ref{1}) be a $2\times 2$ strictly hyperbolic, genuinely nonlinear system, with smooth coefficients.
 Consider a Temple class system satisfying the conclusion of Lemma~\ref{l:41}.  Let $L^*$ be the Lipschitz 
 constant of the semigroup generated by the Temple system.   Then, for any given $\tau>0$, there exists
 $\ve_1>0$ such that the following holds. 

Let $u$ be a front tracking solution of the Cauchy problem (\ref{1})-(\ref{2})  with initial data 
satisfying (\ref{Linfty}). Assume that at time $t=\tau$ a front of size $\sigma_0$ is shifted at rate $\xi_0$.
Then the strengths and shifts of all wave-fronts at any  time $T \in [\tau, 1]$ satisfy
\bel{totsh} \sum_\alpha \bigl|\sigma_\alpha(T)\xi_\alpha(T)\bigr|~\leq~2L^* |\sigma_0\xi_0|\,.
\eeq
\end{lemma}

{\bf Proof.} {\bf 1.}
Without loss of generality, we can assume that $x_0=0$, all wave speeds lie inside the interval $[-1,1]$,  
and the initial datum $\bar u$ is constant
outside the interval $[-2,2]$, say 
\bel{baru0}\bar u(x)~=~\left\{ \bega{rl} u_0^-\quad &\hbox{if}~~x\leq -2,\\[1mm]
u_0^+\quad &\hbox{if}~~x\geq 2.\enda\right.\eeq
Let $x_\alpha(t)$, $\alpha\in \J(t)$ be the location of the jumps in $u(t,\cdot)$.
Call $\sigma_\alpha$ the size of the jump and $i(\alpha)\in \{1,2\}$ its family.
In terms of a system of Riemann coordinates $(w_1, w_2)$, calling $(w_1^{\alpha-}, w_2^{\alpha^-})$
and $(w_1^{\alpha+}, w_2^{\alpha^+})$ the coordinates of the left and right states, we thus define
$$\sigma_\alpha~\doteq~\left\{\bega{rl} w_1^{\alpha+}-w_1^{\alpha-}\quad\hbox{if}~~i(\alpha)=1,\\[1mm]
w_2^{\alpha+}-w_2^{\alpha-}\quad\hbox{if}~~i(\alpha)=2.\enda\right.$$
Here $(w_{1}^- ,w_{2}^- )$ are the Riemann coordinates of $\bar u_0^-$.
\v
{\bf 2.} Consider a front tracking approximate solution to the original system (\ref{1}).  Assume that at the initial time $\tau$ one of the fronts, say located at $x_0$ with size $\sigma_0$, is shifted at rate $\xi_0$.  
At any time $t\in [\tau, T]$,
denote by $\bigl(\sigma_\alpha(t), \xi_\alpha(t)\bigr)$ the corresponding sizes and shifts of the wave-fronts.
For $t\in [\tau,T]$,
we already know that these wave strengths and shifts satisfy the rough bound (\ref{shb3}).
To achieve a sharper bound, we proceed as follows.

Fix a time $s\in [\tau,T]$.  To each front in the above front tracking solution, in addition to  the couple 
$\bigl(\sigma_\alpha(t), \xi_\alpha(t)\bigr)$ we associate another couple 
$\bigl(\Tilde \sigma_\alpha(t), \Tilde \xi_\alpha(t)\bigr)$, inductively defined by the following rules.
\begi
\item For $t\in [\tau, s]$, we set $\bigl(\Tilde \sigma_\alpha(t), \Tilde \xi_\alpha(t)\bigr)=
\bigl(\sigma_\alpha(t),\xi_\alpha(t)\bigr)$.
\item During the remaining  interval  $[s,T]$, the new strengths and shifts $(\Tilde \sigma_\alpha, \Tilde \xi_\alpha)$
are calculated assuming that, at each interaction point, 
the Riemann solver for the Temple class system (\ref{TCS}) constructed in Lemma~\ref{l:41}) is used.
\endi
More precisely,
for any $t\in [s,T]$, call  $\Tilde u(t,\cdot)$ the piecewise constant function
whose Riemann coordinates $(\Tilde  w_1, \Tilde w_2)(t,x)$ are given by
\bel{Te4}\Tilde w_1(t,x)~=~w_1^-+\sum_{x_\alpha(t)< x, ~i(\alpha)=1} \Tilde \sigma_\alpha(t),\qquad\qquad 
\Tilde w_2(t,x)~=~w_2^-+\sum_{x_\alpha(t)< x,~ i(\alpha)=2} \Tilde \sigma_\alpha(t)\,.\eeq
Here $(w_{1}^- ,w_{2}^- )$ are the Riemann coordinates of $\bar u_0^-$ in (\ref{baru0}).
Let $t^*$ be a time where two fronts interact, say at a point $\bar x$.
Call $\tilde u^l, \tilde u^m, \tilde u^r$ the left, middle and right states before the interaction.
By the definition of $\Tilde u$, the Riemann coordinates of these states are computed as in (\ref{Te4}), at a time
$t<t^*$ immediately before the interaction.

Next, call $\Tilde \sigma',\Tilde \sigma''$ the sizes and $\Tilde \xi',\Tilde \xi''$ the
shifts of the incoming fronts.  
By solving the corresponding Riemann problem for the Temple system (\ref{TCS}), we uniquely determine
the sizes  $\Tilde \sigma_1,\Tilde \sigma_2$ and the shifts $\Tilde \xi_1,\Tilde \xi_2$ of the outgoing fronts.

By induction over all interaction times $t^*\in [s,T]$, in a finite number of steps all couples $(\Tilde \sigma_\alpha,
\Tilde\xi_\alpha)$ are determined.

\v
{\bf 3.} In the case where $s=\tau$, the wave strengths and shifts $(\Tilde \sigma_\alpha, \Tilde \xi_\alpha)$
are computed according to the Temple system, at every interaction point. 
We claim that  in this setting the analysis in Section~\ref{sec:4} remains valid. Namely, for every $T\in [\tau,1]$ these strengths and shifts 
satisfy
\bel{diff2}\sum_{\alpha\in \J(T)} \bigl| \Tilde \sigma_\alpha(T)\Tilde\xi_\alpha(T)\bigr|~\leq~
L^* |\sigma_0\xi_0|,\eeq
up to an arbitrarily small increase of the constant $L^*$ depending on $\ve_1$.

Indeed the shifts $\xi_\alpha(T)$ can again be computed by the explicit formulas (\ref{shift1}), (\ref{shift2}), (\ref{shift3a}) and (\ref{shift3b}), where $u^\pm$ and $u_0^\pm$ are replaced by $\Tilde u^\pm$, $\Tilde u^\pm_0$. 
This yields the estimates \eqref{sb1}, \eqref{sb2}, \eqref{sb3a} and \eqref{sb3b}. 

In order to conclude, it remains to check that we have a uniform bound of the form
\bel{e:tv-mix}
\sum_{x_\alpha\in  [x_0-T,\, x_0+T]}  |\Tilde \sigma_\alpha|~\leq~\Tilde C_0.
\eeq
Indeed, by \eqref{e:tv-gnl} we can estimate
\bel{e:diff-str}
\begin{split}
\sum_{x_\alpha\in  [x_0-T,\, x_0+T]}  |\Tilde \sigma_\alpha|~\leq &~ \sum_{x_\alpha\in  [x_0-T,\, x_0+T]}  | \sigma_\alpha| + \sum_{x_\alpha\in  [x_0-T,\, x_0+T]}  |\Tilde \sigma_\alpha - \sigma_\alpha|\\
\leq & ~ C_0 +  \sum_{x_\alpha\in  [x_0-T,\, x_0+T]}  |\Tilde \sigma_\alpha - \sigma_\alpha|.
\end{split}
\eeq
In order to estimate the last sum, consider the interaction functional $Q(t)=\sum_{(\alpha,\beta) \in \mathcal A}|\sigma_\alpha \sigma_\beta|$, where $\mathcal A$ denotes the set of all approaching waves with $x_\alpha(t), x_\beta(t) \in [x_0-2T+t,\, x_0+2T-t]$.
Thanks to the bounds (\ref{ies1}), (\ref{ies3}) and the uniform bound $|\sigma| = \O(1)\cdot \sqrt{\ve_1}$ 
on the strength of each  front, there is a constant $C>0$ such that
\[
\Gamma(t)~ \doteq \sum_{x_\alpha\in  [x_0-2T+t,\, x_0+2T-t]} |\Tilde \sigma_\alpha(t) - \sigma_\alpha(t)| + C \sqrt{\ve_1} Q(t)
\]
is decreasing. In particular 
\[
 \sum_{x_\alpha\in  [x_0-T,\, x_0+T]}  |\Tilde \sigma_\alpha - \sigma_\alpha| ~\leq~ \Gamma(0) ~=~ \O(1) 
 \cdot \sqrt{\ve_1} \tau^{-2}.
\]
Inserting the above estimate in \eqref{e:diff-str}, we obtain \eqref{e:tv-mix} and hence \eqref{diff2}.
\v
{\bf 4.}
To complete the proof of Lemma~\ref{l:61} we will show that, at the terminal time $T\in [\tau,1]$, 
if $\ve_1>0$ was chosen small enough,
the difference between the two sets of couples (computing the sizes and shifts of all wavefronts in two different ways)
is bounded by 
\bel{diff1}\sum_{\alpha\in \J(T)} \bigl| \sigma_\alpha(T) \xi_\alpha(T)  - \Tilde \sigma_\alpha(T) \Tilde\xi_\alpha(T) \bigr|~\leq~
L^* |\sigma_0\xi_0|.\eeq

\begin{figure}[ht]
  \centerline{\hbox{\includegraphics[width=10cm]{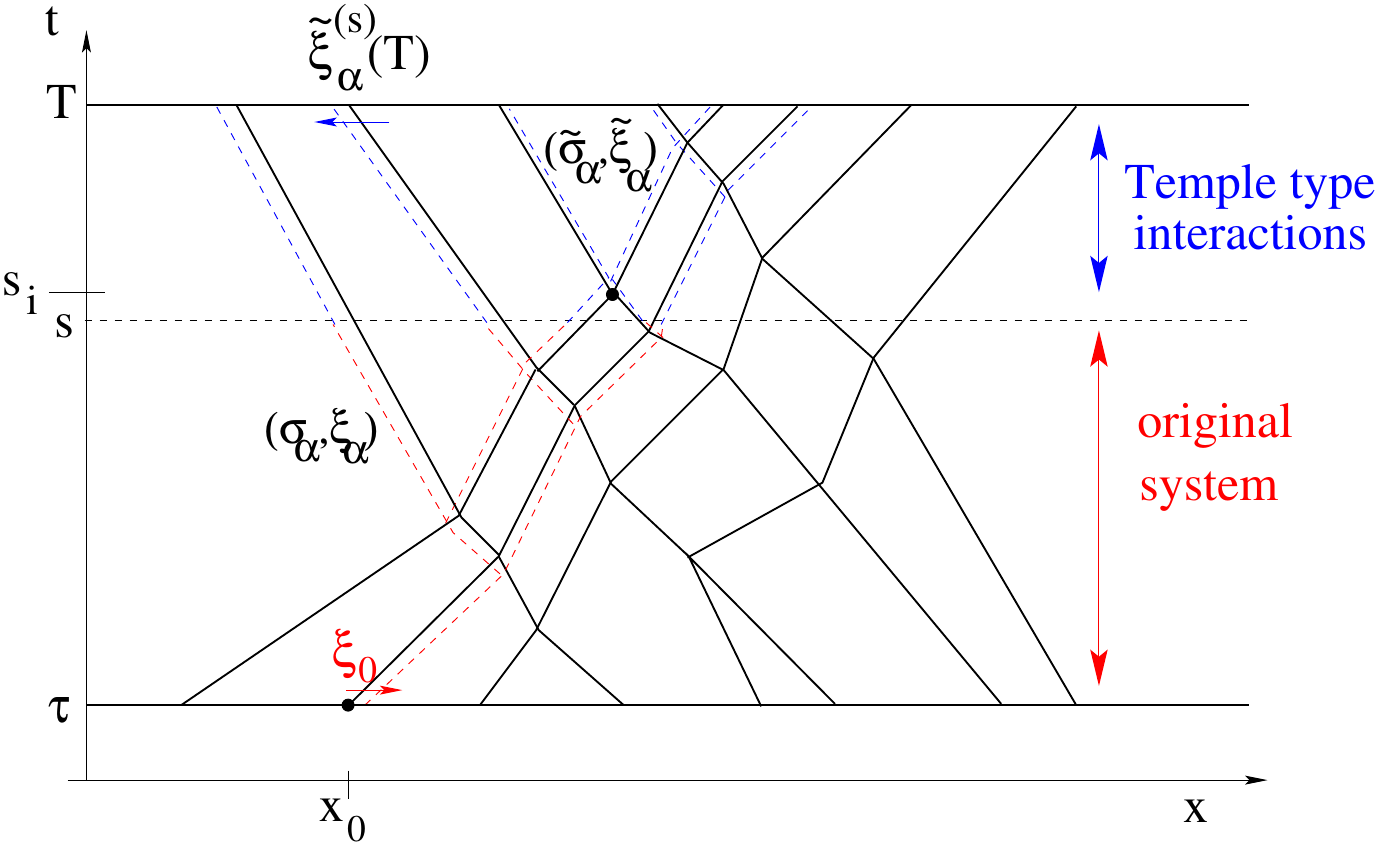}}}
  \caption{\small The homotopy argument yielding the inequality (\ref{diff1}). 
  For a given  intermediate time $s\in [\tau,T]$    we calculate the strengths and the shifts of 
  the fronts according to the original system of conservation laws (\ref{1}) for interactions at times
  $\tau<s_i<s$, and according to the Temple system (\ref{TCS}) at interaction times $s<s_i<T$.}\label{f:ag88}
\end{figure}

 The inequality (\ref{diff1}) will be proved by a homotopy argument.
For every $s\in [\tau, T]$,  
call $\Big(\Tilde \sigma_\alpha^{(s)} (T),\Tilde \xi_\alpha^{(s)}(T)\Big)$ the strengths and shifts of the wavefronts at the terminal time
$T$, computed as in  step  {\bf 2}.    Notice that here the superscript keeps track of the dependence 
on the intermediate time $s\in [\tau,T]$.
These values can change only across the finitely many 
 times $s_i\in[\tau,T]$ when
two wave-fronts in the original front tracking solution $u$ interact.

We denote by $\I_d$ the set of indexes $i$ such that an interaction of fronts of different families happens at time $s_i$ and by $\I_s$ the set of indexes $i$ such that an interaction of fronts of the same family happens at time $s_i$.
 Call $\sigma_i',\sigma_i''$ and $\xi'_i,\xi_i''$ respectively the sizes and shifts of the two approaching fronts before  the interaction at time $s_i$.  Two cases must be considered.
 
CASE  1: $i\in \I_d$. 
We claim that
\bel{Tshd} \sum_{\alpha} \Big| \Tilde \sigma_\alpha^{(s_i+)}(T)\Tilde  \xi_\alpha^{(s_i+)}(T) \Big| -  \sum_{\alpha}\Big|\Tilde \sigma_\alpha^{(s_i-)}(T) \Tilde \xi_\alpha^{(s_i-)}(T) \Big|~\leq~C_2\sqrt{\ve_1}\, 
\bigl( |\xi_i'|+|\xi_i''|\bigr) |\sigma_i'\sigma_i''|.
\eeq
Indeed, the strengths $\sigma_\alpha^{(s_i\pm)}$ and the shifts  $\xi_\alpha^{(s_i\pm)}$ 
are computed for the same Temple class system on the interval $t\in [s_i+\delta, T]$  (with $\delta>0$ small), but with 
different initial data at time $t=s_i+\delta$.   
In the first case, the strengths and shifts of fronts emerging from the
interaction at time $s_i$ are computed according to the system (\ref{1}), while in the 
second case they are computed according to the Temple system  (\ref{TCS}).  
Using Lemma~\ref{l:52} with $\ve =\O(1)\cdot \sqrt{\ve_1}$, the differences in the strengths and shifts
at time $s_i+\delta$, immediately after the interaction, satisfy the bounds (\ref{bT1})-(\ref{bT2}), with $\ve=\O(1)\cdot\sqrt{\ve_1}$.

We can now use the explicit formulas (\ref{shift1}), (\ref{shift2}), (\ref{shift3a}) and (\ref{shift3b}) to compute
the corresponding shifts at the terminal time $t=T$.  Notice that these are linear w.r.t.~the initial shifts, and
Lipschitz continuous w.r.t.~the initial and terminal states $u_0^\pm$ and $u^\pm$.
Since the total strength of all fronts at time $T$ is uniformly bounded, this yields the estimate (\ref{Tshd}).
\v
CASE 2: $i\in \I_s$.   In this case, by Lemma~\ref{l:52} we have the bounds (\ref{bT3})--(\ref{bT5}).
The same arguments as in the previous case now yield only the weaker estimate
\bel{Tshs} \sum_{\alpha} \Big|\Tilde  \sigma_\alpha^{(s_i+)}(T) \Tilde \xi_\alpha^{(s_i+)}(T) \Big| -  \sum_{\alpha}\Big|\Tilde \sigma_\alpha^{(s_i-)}(T) \Tilde \xi_\alpha^{(s_i-)}(T) \Big|~\leq~C_2\, 
\bigl( |\xi_i'|+|\xi_i''|\bigr) |\sigma_i'\sigma_i''|.
\eeq

By (\ref{Qe}), in the front tracking solution for the original system (\ref{1}), the sum of all interactions with shifts is uniformly bounded.
In particular 
\[
\sum_{i \in \I_d} \bigl( |\xi_i'|+|\xi_i''|\bigr) |\sigma_i'\sigma_i''| ~\leq~\mathcal{Q}^\xi~\leq~C_3 |\sigma_0 \xi_0|
\]
for some constant $C_3$. 
On the other hand, 
by Lemma~\ref{l:32} the total amount of interactions with shifts among fronts of the same family  is bounded by
\[
\sum_{i \in I_s} \bigl( |\xi_i'|+|\xi_i''|\bigr) |\sigma_i'\sigma_i''| ~\leq~\mathcal{Q}_{same}^\xi~\leq~C_4 \sqrt{\ve_1}\, |\sigma_0 \xi_0|.
\]
 We now observe that, when $s=T$, by definition one trivially has $\sigma_\alpha(T)= \Tilde \sigma_\alpha^{(T)}(T)$ and 
 $\xi_\alpha(T)=\Tilde \xi_\alpha^{(T)}(T) $ for every $\alpha\in \J(T)$.
By choosing $\ve_1>0$ small enough, we can now achieve the bound
\bel{Tsh3}\bega{l} \ds
 \sum_\alpha \Big| \Tilde \sigma_\alpha^{(T)}(T)\Tilde  \xi_\alpha^{(T)}(T) \Big| -  \sum_\alpha\Big|\Tilde \sigma_\alpha^{(\tau)}(T) \Tilde \xi_\alpha^{(\tau)}(T) \Big|\\[4mm]
\qquad \ds\leq~\bigg(\sum_{i\in \I_d}+ \sum_{i\in \I_s}\bigg)
 \bigg(\sum_\alpha \Big|\Tilde  \sigma_\alpha^{(s_i+)}(T)\Tilde  \xi_\alpha^{(s_i+)}(T) \Big| -  \sum_\alpha\Big|\Tilde \sigma_\alpha^{(s_i-)}(T) \Tilde \xi_\alpha^{(s_i-)}(T) \Big|\bigg)\\[6mm]
 \ds\qquad \leq~ C_2\sqrt{\ve_1} C_3 |\sigma_0 \xi_0| + C_2 C_4 \sqrt{\ve_1} |\sigma_0 \xi_0|~\leq~L^* |\sigma_0 \xi_0|.\enda
\eeq
Combining (\ref{Tsh3}) with (\ref{diff2}) the proof is completed.
\endproof

\section{Proof of Theorem~\ref{t:11}}
\label{sec:7}
\setcounter{equation}{0}
{\bf 1.} Let (\ref{TCS}) be a Temple class system satisfying the conclusion of Lemma~\ref{l:41}, generating a
semigroup with Lipschitz constant $L^*$.  
Let $q>1$ be a constant large enough so that 
\bel{loga}  {\ln (2L^*)\over\ln q}~<~\alpha\,.\eeq
Taking $\tau=1/q$, according to Lemma~\ref{l:61}, choosing $\ve_1>0$
small enough we ensure that the 
Lipschitz constant of the semigroup $S$ generated by (\ref{1}) on the domain $\D^{(1/q)}$ 
is $\leq 2L^*$.

\begin{figure}[ht]
  \centerline{\hbox{\includegraphics[width=9cm]{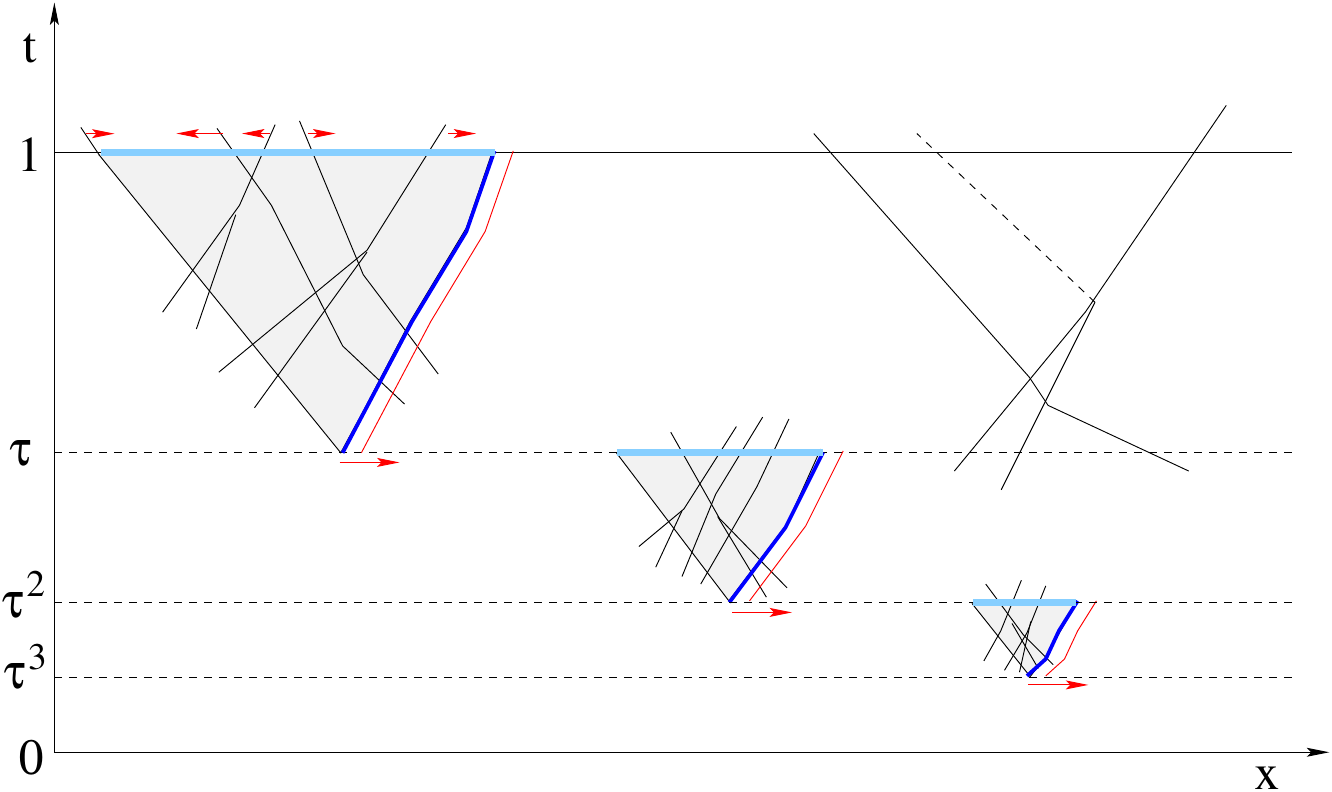}}}
  \caption{\small The amount of waves affected by a shifted front during the time interval
  $[\tau,1]$ is of the same order of magnitude as the amount of waves affected by a shifted front during the time interval
  $[\tau^2,\tau]$.    Same for $[\tau^3,\tau^2]$, etc...}\label{f:ag65}
\end{figure}

\v
{\bf 2.} 
By rescaling, we obtain that the same  Lipschitz constant is valid 
on any interval $[q^{-k}, q^{-k+1}]$, namely
\bel{Lip5}
\bigl\| S_t\bar u-S_t\bar v\bigr\|_{\L^1}~\leq~2L^*\cdot \|\bar u - \bar v\|_{\L^1}\,,\qquad\qquad 
\forall \bar u,\bar v\in \D^{(q^{-k})}, ~~t\in [q^{-k}, q^{-k+1} ].\eeq
In turn, an inductive argument yields
\bel{Lip6}
\bigl\| S_t\bar u-S_t\bar v\bigr\|_{\L^1}~\leq~(2L^*)^k\cdot \|\bar u - \bar v\|_{\L^1}\,,\qquad\qquad 
\forall \bar u,\bar v\in \D^{(q^{-k})}, ~~t\in [q^{-k}, 1].\eeq
In other words,
$$L(\tau)~\leq~(2L^*)^k\qquad\qquad\hbox{for}\quad \tau\geq q^{-k}.$$
Since this holds with 
$$k~\leq~1 + \ln_q {1\over\tau} ~=~1- {\ln\tau\over\ln q}\,,$$
we conclude
\bel{Ltau}L(\tau)~\leq~\exp\left\{ \left( 1-{ \ln \tau\over\ln q}\right) \cdot \ln (2L^*)\right\}
~\leq~2L^*\cdot \tau^{-\beta},\qquad\qquad \beta \doteq {\ln (2L^*)\over\ln q}~<~\alpha\,.\eeq
\v
{\bf 3.} To prove the H\"older continuous dependence on initial data, 
consider two initial data $\bar u_1,\bar u_2$.   Then for any $0<\tau<t\leq 1$, by (\ref{dista})  and (\ref{Ltau}) 
it follows
$$\bigl\|u_1(t)-u_2(t)\bigr\|_{\L^1}
~\leq~L(\tau) \cdot \bigl( \|\bar u_1-\bar u_2\|_{\L^1} + 2C \tau^\alpha\bigr).$$
If $t\geq \|\bar u_1-\bar u_2\|_{\L^1}^{1/\alpha}$, choosing 
$\tau\doteq \|\bar u_1-\bar u_2\|_{\L^1}^{1/\alpha}\bigr\}$ we achieve
\bel{u12d}
\bigl\|u_1(t)-u_2(t)\bigr\|_{\L^1}~\leq~2L^* \|\bar u_1-\bar u_2\|_{\L^1}^{-\beta/\alpha} \cdot  \|\bar u_1-\bar u_2\|_{\L^1}
\cdot (1+ 2C).
\eeq
This yields H\"older continuity, with exponent $\gamma = 1-(\beta/\alpha)$.
On the other hand, if $t\leq \|\bar u_1-\bar u_2\|_{\L^1}^{1/\alpha}$, we have the Lipschitz bound
\bel{u12s}
\bigl\|u_1(t)-u_2(t)\bigr\|_{\L^1}~
\leq~ \|\bar u_1-\bar u_2\|_{\L^1} + 2Ct^\alpha~
\leq~ \|\bar u_1-\bar u_2\|_{\L^1} \cdot (1+2C).
\eeq
Note that, by choosing $q$ suitably large, the constant $\beta$ in (\ref{Ltau}) can be rendered arbitrarily small.
In turn, the H\"older exponent $\gamma = 1-(\beta/\alpha)$ becomes arbitrarily close to 1.
\v
{\bf 4.} Since both characteristic fields are assumed to be genuinely nonlinear, by  the recent
results in
\cite{BDL, BGu} the Lax condition suffices to conclude the uniqueness of solutions
that satisfy the decay estimate (\ref{decal}).
\endproof

\section{A domain of initial data with fast decay}
\label{sec:8}
\setcounter{equation}{0}
Aim of this section is to prove that the condition {\bf ($\Tilde P_\alpha$)} 
introduced in Section~1 yields fast decay of the total variation.
This provides a domain of initial data, possibly with unbounded variation,
where the hyperbolic system (\ref{1}) has unique solutions. As before, w.l.o.g.~we assume that all
wave speeds are contained in the interval $[-1,1]$.

The proof of Theorem \ref{t:12} relies on the following proposition.
\v
\begin{prop}\label{p:BCM}
There exist constants $C,K,\varepsilon_1>0$ such that, if $\|\bar u\|_{\bf L^\infty}< \varepsilon_1$, then \eqref{1} has a solution $u$ satisfying the following estimates, for every $t > 0$ and every $a,b \in \R$:
\begin{equation}\label{e:TV1}
\TV \bigl\{ u(t);[a,b]\bigr\}~\le~ C\left(1+\frac{b-a}{t}\right).
\end{equation}
If moreover $\TV\{ \bar u; ~[a-  t, b +   t]\} \le K$ for some $t>0$ and $a,b\in \R$, then
\begin{equation}\label{e:TV2}
\TV\bigl\{ u(t);~[a,b] \bigr\} ~\le~ C\, \TV\bigl\{ \bar u; ~[a-  t, b +   t ]\bigr\}.
\end{equation}
\end{prop}

Indeed, 
The first part of the statement, and in particular \eqref{e:TV1}, is contained in \cite{BCM}. The estimate \eqref{e:TV2} follows from the analogous estimates for front tracking approximate solutions as in \cite{BC95}. 

\v
{\bf Proof of Theorem \ref{t:12}.}

{\bf 1.} Fix $t\in \,]0,1]$ and, for every $k \in \mathbb{Z}$ and $\tau \in [0,t]$, define the interval
\[
I_k(\tau)\,=\,\bigl[kt -   (t-\tau), (k+1)t +   (t-\tau)\bigr].
\]
We denote by $\J_{[a,b]}\subset \mathbb{Z}$ the set of indices $k$ for which $I_k(t)\cap [a , b ]\ne \emptyset$ and by $\J\subset \J_{[a,b]}$ the set of indices $k$ for which $\TV\{\bar u; I_k(0)\}>K$.

By assumption, the initial data $\bar u$ satisfies the condition ($\widetilde {\bf P}_\alpha$).
Setting $\lambda=t$, we thus call $V^t=V^\lambda$ the set  which satisfies the conditions (\ref{mop2})--(\ref{Nb}),
in connection with the interval $I=  [a - t, b +  t] \supset \bigcup_{k \in \J_{[a,b]}}I_k(0)$.  
Finally,  we denote by $\J'\subset \J_{[a,b]}$ the set of indices $k$ for which $I_k(0)\cap V^\lambda \ne 0$.
Since the family $\{I_k(t)\}_{k\in \J_{[a,b]}}$ covers $[a,b]$, we have 
\bel{TV9}
\TV \bigl\{u(t);\,[a,b]\bigr\} \le \sum_{k \in \J\cup \J'}\TV \bigl\{u(t);\, I_k(t)\bigr\} 
+ \sum_{k \in \J_{[a,b]} \setminus (\J\cup \J')}\TV \bigl\{u(t);\,I_k(t)\bigr\}.
\eeq
The above  two sums will be estimated separately.

By \eqref{e:TV1}, for each $k \in \mathbb{Z}$ one has
\begin{equation}\label{e:TVIk}
\TV \bigl\{ u(t); \,I_k(t) \bigr\} ~\le ~2C.
\end{equation}
We now estimate the cardinality of $\J \cup \J'$.  Since every $x  \in \R$ satisfies the inclusion $x \in I_k(0)$ for at most $2+2 $ indices $k$ and the distance between two disjoint intervals $I_k(0),I_{k'}(0)$ is at least $t$, each connected component $]a_i,b_i[$ of $V^\lambda$ intersects at most $4 (1+ \frac{b_i-a_i}t)$ intervals $I_k(0)$. Summing over $i=1,\ldots,N$, the assumptions \eqref{mop2} and  \eqref{Nb}  imply that $\J'$ contains at most $8\Tilde Ct^{\alpha -1}$ indices.

Using again that every $x  \in \R$ satisfies the inclusion $x \in I_k(0)$ for at most $2+2 $ indices $k$,   and moreover
$$\TV\bigl\{\bar u;~ [a-  t,b+  t] \setminus V^\lambda \bigr\}~ \le~ \tilde C t^{\alpha -1},$$ 
we conclude that
\begin{equation}\label{e:TVJ'c}
\sum_{k \in \J_{[a,b]} \setminus \J'} \TV \bigl\{u(t); \,I_k(0)\bigr\} ~\le ~ \left(2+ 2 \right)\Tilde C t^{\alpha -1}.
\end{equation}
Therefore the cardinality of $\J\setminus \J'$ is bounded by $ \frac{2+ 2 }{K}\Tilde C t^{\alpha -1}$.
In particular there is $C_{K}>0$ such that the cardinality of $\J$ is bounded by $C_K \Tilde C t^{\alpha -1}$.
Therefore, by \eqref{e:TVIk} one has
\begin{equation}\label{e:JJ'}
 \sum_{k \in \J\cup \J'}\TV \bigl\{u(t);\, I_k(t)\bigr\}~\le~ 2  C_KC \Tilde C t^{\alpha -1}.
\end{equation}
Finally, by \eqref{e:TV2} and \eqref{e:TVJ'c}, we have
\begin{equation}\label{e:JJ'c}
\sum_{k \in \J_{[a,b]} \setminus (\J\cup \J')}\TV \bigl\{u(t); I_k(t)\bigr\}~\le~ 
C \sum_{k \in \J_{[a,b]}\setminus (\J\cup \J')}\TV \bigl\{\bar u;\, I_k(0)\bigr\}~\le~  \left(2+ 2 \right)C \Tilde C t^{\alpha -1}.
\end{equation}
Using the two bounds \eqref{e:JJ'} and \eqref{e:JJ'c} in (\ref{TV9}), we conclude that the total variation decays like $t^{\alpha-1}$, proving the theorem.
\endproof

{ \bf Acknowledgements. } The research of A.~Bressan  was partially supported by NSF  with
grant  DMS-2306926 ``Regularity and approximation of solutions to conservation laws". 
E.~Marconi is a member of GNAMPA of the “Istituto Nazionale di Alta Matematica F. Severi” and is  partially supported by H2020-MSCA-IF “A Lagrangian approach: from conservation laws to line-energy Ginzburg-Landau models”.

\v

\end{document}